\documentclass[oneside]{amsart}
\usepackage[parfill]{parskip}
\usepackage{amsmath, latexsym, amsthm, mathtools,subfigure,float}
\usepackage[title]{appendix}
\usepackage[sort,nocompress]{cite}

\usepackage[margin=1.15in]{geometry}

\providecommand{\N}{\mathbb{N}}
\providecommand{\R}{\mathbb{R}}

\providecommand{\Z}{\mathbb{Z}}

\DeclarePairedDelimiter\floor{\lfloor}{\rfloor}

\providecommand{\abs}[1]{\left\vert#1\right\vert}

\providecommand{\set}[1]{\left\{#1\right\}}
\providecommand{\paren}[1]{\left( #1 \right)}

\providecommand{\brac}[1]{\left [ #1 \right]}
\providecommand{\set}[1]{\left { #1 \right}}

\newcommand{\PH}{\mathit{PH}}

\newtheorem{Theorem}{Theorem}

\newtheorem{Definition}[Theorem]{Definition}

\begin{document}
\title[Fractal Dimension Estimation with Persistent Homology]{Fractal Dimension Estimation with Persistent Homology: \\ A Comparative Study}


\author{Jonathan Jaquette}
\address[J. Jaquette]{Department of Mathematics, 
	Brandeis University
	Waltham, MA 02453}
\email{jjaquette@brandeis.edu}
\thanks{Research of the first author was supported in part by NIH T32 NS007292 and NSF DMS-1440140 while the author was in residence at the Mathematical Sciences Research Institute in Berkeley, California, during the Fall 2018 semester.}

\author{Benjamin Schweinhart}
\thanks{Research of the second author was supported in part by a NSF Mathematical Sciences Postdoctoral Research Fellowship under award number DMS-1606259.} 
\address[B. Schweinhart]{Department of Mathematics,  Ohio State University Columbus, OH 43210}
\email{schweinhart.2@osu.edu}

\date{August 2019}

\begin{abstract}
 We propose that the recently defined persistent homology dimensions are a practical tool for fractal dimension estimation of point samples. We implement an algorithm to estimate the persistent homology dimension, and compare its performance to classical methods to compute the correlation and box-counting dimensions in examples of self-similar fractals, chaotic attractors, and an empirical dataset. The performance of the $0$-dimensional persistent homology dimension is comparable to that of the correlation dimension, and better than box-counting.  

\end{abstract}

\maketitle

\section{Introduction}
Loosely speaking, fractal dimension measures how local properties of a set depend on the scale at which they are measured. The Hausdorff dimension was perhaps the first precisely defined notion of fractal dimension~\cite{1918hausdorff,2003edgar}. It is difficult to estimate in practice, but several other more computationally practicable definitions have been proposed, including the box-counting~\cite{1928bouligand} and correlation~\cite{1983grassberger} dimensions. These notions are in-equivalent in general.

Popularized by Mandelbrot~\cite{1977mandelbrot,1982mandelbrot}, fractal dimension has applications in a variety of fields including  materials science~\cite{1999davies,2008yu,2003hu}, biology~\cite{2000baish,2009lopes,2013peng}, soil morphology~\cite{1991rieu}, and the analysis of large data sets~\cite{2000barbara,2010traina}. It is also important in pure mathematics and mathematical physics, in disciplines ranging from dynamics~\cite{1980takens} to probability~\cite{2008beffara}. In some applications, it is necessary to estimate the dimension of a set from a point sample. These include earthquake hypocenters and epicenters~\cite{2007kagan,1998harte}, rain droplets~\cite{1990lovejoy,2000gabella}, galaxy locations~\cite{1991guzzo}, and chaotic attractors~\cite{1980takens,1983grassberger}.

\begin{figure}
\centering  
\subfigure[]{
\label{fig:sierpinski}
\includegraphics[height=0.27\linewidth]{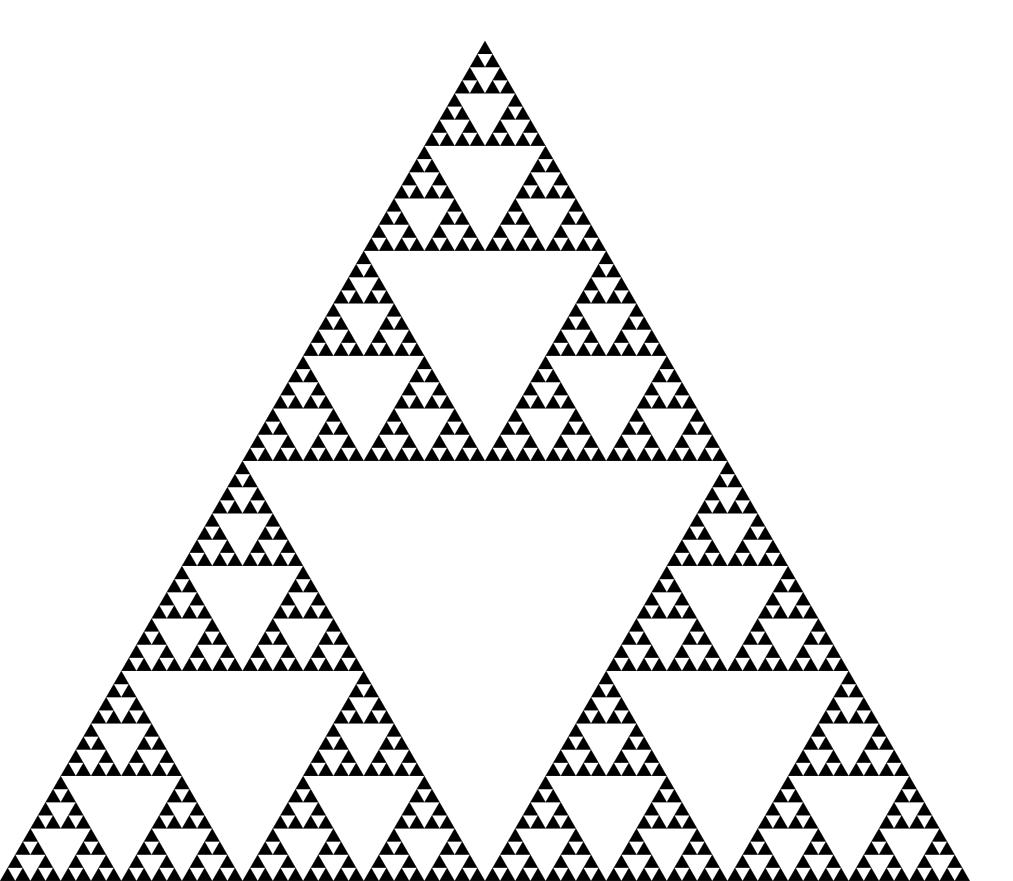}
}
\hspace{.025\linewidth}
\subfigure[]{
\label{fig:cantordust}
\includegraphics[height=0.27\linewidth]{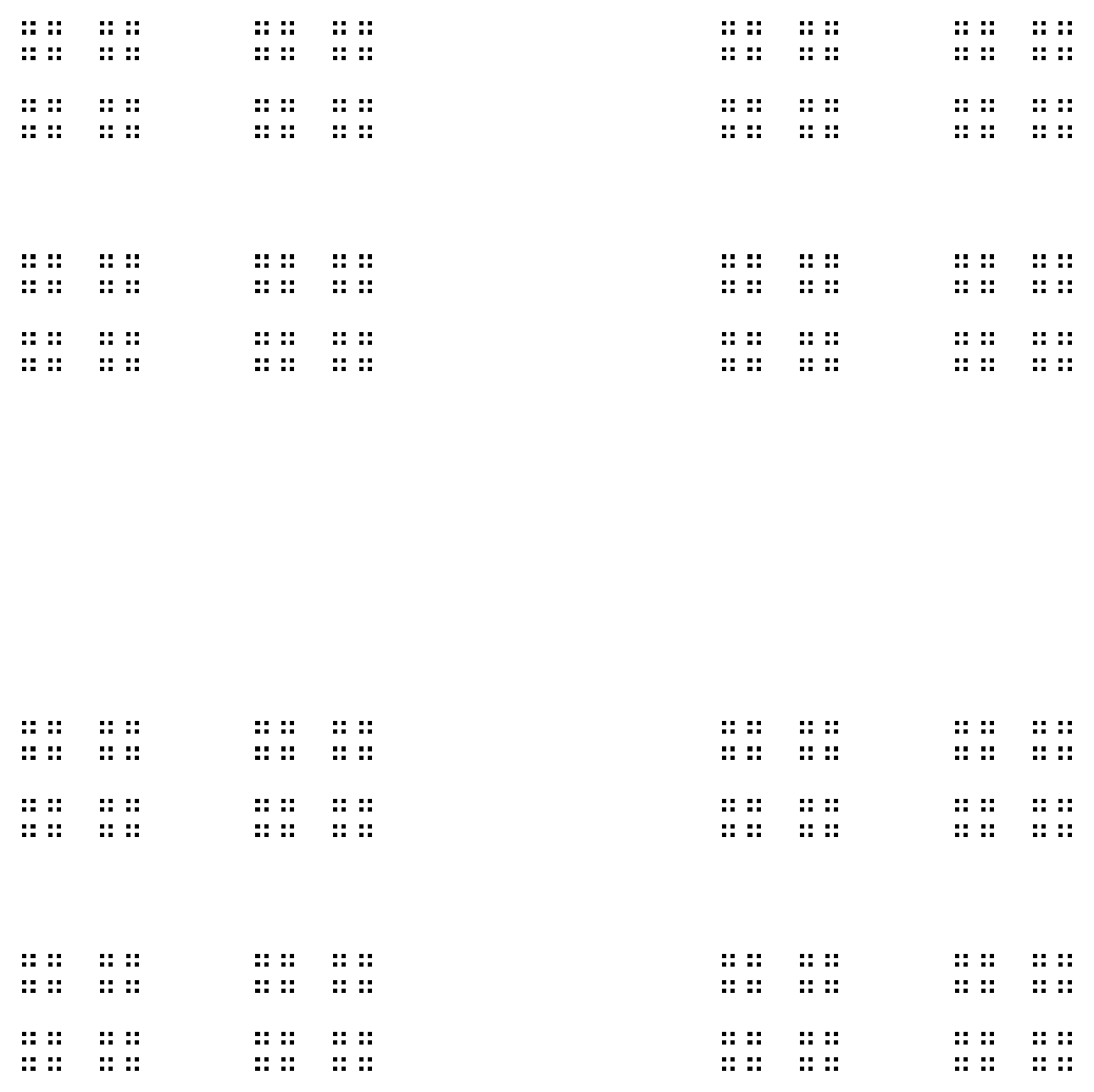}}
\hspace{.025\linewidth}
\subfigure[]{
\label{fig:CXI}
\includegraphics[height=0.27\linewidth]{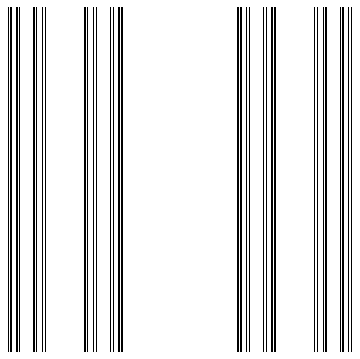}}

\vspace{1\baselineskip}

\subfigure[]{
	\label{fig:Henon}
	\includegraphics[height=0.27\linewidth]{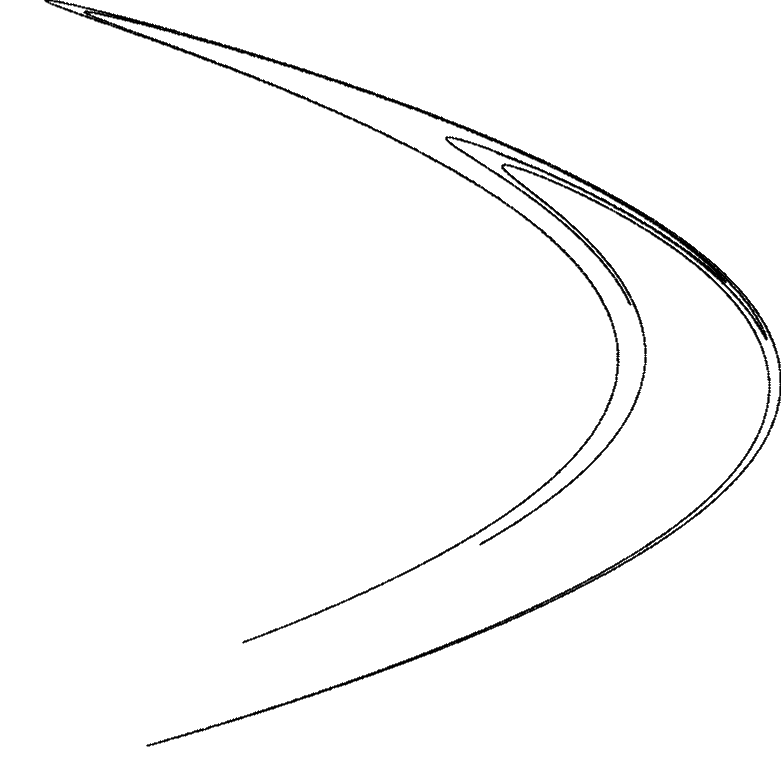}
}
\hspace{.025\linewidth}
\subfigure[]{
	\label{fig:Ikeda}
	\includegraphics[height=0.27\linewidth]{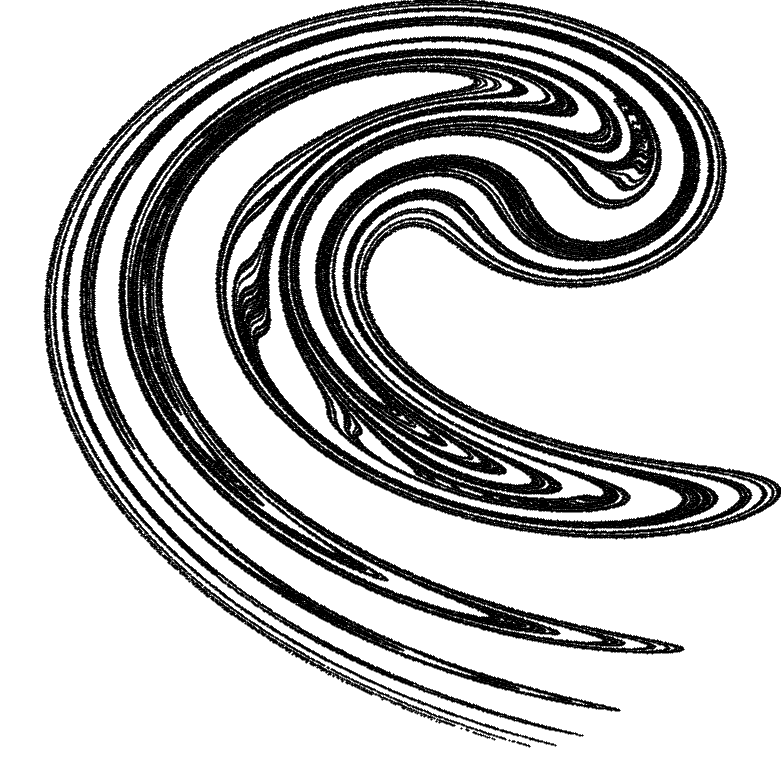}}
\hspace{.025\linewidth}
\subfigure[]{
	\label{fig:Rulkov}
	\includegraphics[height=0.27\linewidth]{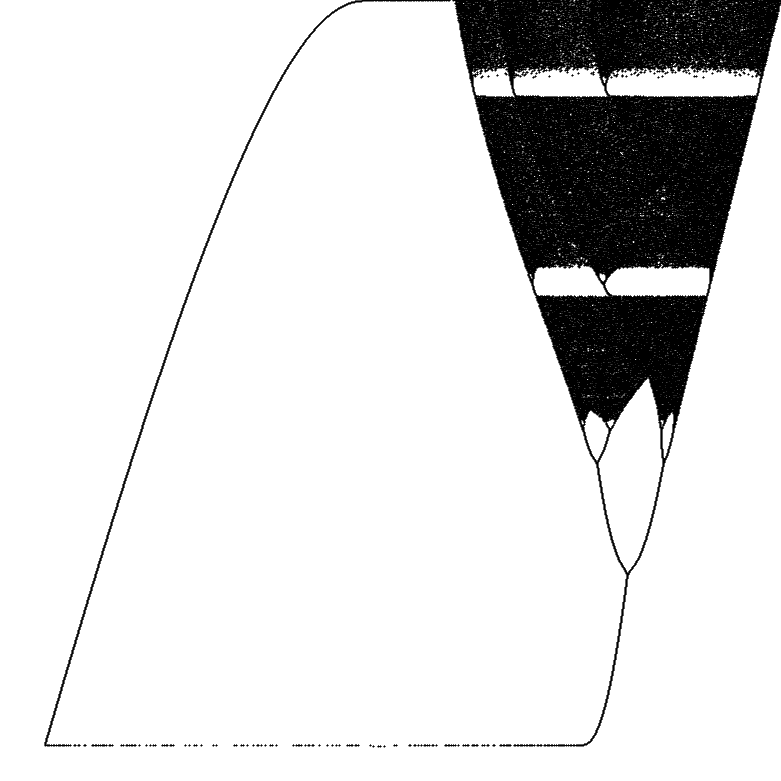}}

		\caption{\label{fig:fractals}Three self-similar fractals: (a) the Sierpinski Triangle, (b) the Cantor dust, and (c) the Cantor set cross an interval; and three chaotic attractors arising from the: (d) H\'enon map (e) Ikeda map and (f) Rulkov map.}
\end{figure}

We propose that the recently defined persistent homology dimensions~\cite{2019adams,2019schweinhart} (Definition~\ref{defn:phdim} below) are a practical tool for dimension estimation of point samples.  Persistent homology~\cite{2002edelsbrunner} quantifies the shape of a geometric object in terms of how its topology changes as it is thickened; roughly speaking, $i$-dimensional persistent homology ($\PH_i$) tracks $i$-dimensional holes that form and disappear in this process. Recently, it has found many applications in fields ranging from materials science~\cite{2016hiraoka,2017saadatfar} to biology~\cite{2015giusti,2014xia}. There are further applications of persistent homology in machine learning, for which different methods of vectorizing persistent homology have been defined~\cite{2017adams,2015bubenik}. In most applications, larger geometric features represented by persistent homology are of greatest interest. However, in the current context it is the smaller features --- the ``noise'' --- from which the dimension may be recovered. 

For a finite subset of a metric space $\set{x_1,\ldots,x_n}$ and a positive real number $\alpha$ define the $\alpha$ weight $E_\alpha^i\paren{x_1,\ldots,x_n}$ to be the sum of the lengths of the (finite) $\PH_i$ intervals to the $\alpha$ power:
\[E_\alpha^i\paren{x_1,\ldots,x_n}=\sum_{I\in\PH_i\paren{x_1,\ldots,x_n}}\abs{I}^\alpha\,.\]
If $\mu$ is a measure on bounded metric space, $\set{x_j}_{j\in\mathbb{N}}$ are i.i.d.~samples from $\mu,$ and $i$ is a natural number the $\PH_i$-dimension $\text{dim}_{\PH_i^\alpha} \paren{\mu}$ measures how $E_\alpha^i\paren{x_1,\ldots,x_n}$ scales as $n\rightarrow\infty.$ In particular, if $E_\alpha^i\paren{x_1,\ldots,x_n}\approx n^{\frac{d-\alpha}{d}},$ then $\text{dim}_{\PH_i^\alpha} \paren{\mu}=d.$ See below for a precise definition.

We estimate the persistent homology dimension of several examples, and compare its performance to classical methods to estimate the correlation and box-counting dimensions.  We study the convergence of estimates as the sample size increases, and the variability of the estimate between different samples. The examples we consider are from three broad classes: self-similar fractals, chaotic attractors, and empirical data. The sets in the former class (shown in Figure~\ref{fig:sierpinski}--~\ref{fig:CXI}) have known dimensions, and are regular in a sense that implies that the various notions of fractal dimension agree for them. That is, they have a single well-defined ``fractal dimension.'' This makes it easy to compare the performance of the different dimension estimation techniques. In those cases, our previous theoretical results~\cite{2019schweinhart} imply that at least the zero-dimensional version of persistent homology dimension will converge to the true dimension. 

The second class of fractals we study are strange attractors arising from chaotic dynamical systems, see e.g. Figure \ref{fig:Henon}-\ref{fig:Rulkov}. Generically, these sets are not known to be regular, and the various definitions of fractal dimension may disagree for them. Finally, we apply the dimension estimation techniques to the Hauksson--Shearer Southern California earthquake catalog~\cite{2012hauksson}.
 
We also propose that a notion of persistent homology complexity due to MacPherson and Schweinhart~\cite{2012macpherson} may be a good indicator of how difficult it is to estimate the correlation dimension or persistent homology dimension of a shape.

\pagebreak[3]

In summary, we conclude the following.
\begin{itemize}

\item \textbf{Effectiveness}. Based on our experiments, the $\PH_0$ and correlation dimensions perform comparably well. In cases where the true dimension is known, they approach it at about the same rate. In most cases, the box-counting, $\PH_1,$ and $\PH_2$ dimensions perform worse.
\item \textbf{Efficiency}. Computation of the $\PH_0$ dimension is fast and comparable with the correlation and box-counting dimensions. The $\PH_1$ dimension is reasonably fast subsets of $\mathbb{R}^2,$ but the $\PH_1$ and $\PH_2$ dimensions are quite slow for point clouds in $\mathbb{R}^3,$ and computations for sets with higher ambient dimension are impractical.
\item \textbf{Equivalence}. For a large class of regular fractals the $\PH$ dimension coincides with various classical definitions of  fractal dimension. However, there are fractal sets for which these definitions do not agree, with sometimes surprisingly large differences (e.g. the Rulkov and Mackey-Glass attractors).  Within the class of $\PH_i$ dimensions, there is variation between dimension estimates coming from different choices of the homological dimension or scaling weight. All this said, it can be difficult to determine whether dimensions truly disagree, or if the convergence is very slow.
\item \textbf{Error} Error estimates, whether they are the ``statistical test error'' of the correlation dimension or the empirical standard deviation of estimates between trials, do not contain meaningful information about the difference between the dimension estimate and the true dimension. In general, it is difficult to tell whether a dimension estimate has approached its limiting value.
\item \textbf{Ease-of-use} We found one simple rule for fitting a power law to estimate the $\PH_0$ which worked well for all examples, in contrast to the correlation dimension and (especially) the box-counting dimension. 

\end{itemize}

In the following, we briefly survey previous work comparing different methods for the estimation of fractal dimension (Section~\ref{sec:background}), outline the different methods considered here (Section~\ref{sec:methods}), and discuss the results for each example (Sections~\ref{sec:examples_discrete},~\ref{sec:examples_dynamics}, and~\ref{sec:earthquake}).

\subsection{Background and Previous Work}
 \label{sec:background}

Many previous studies have compared different methods for fractal dimension estimation. These include surveys focusing on applications to chaotic attractors~\cite{1990theiler}, medical image analysis~\cite{2009lopes}, networks~\cite{2018rosenberg}, and time series and spatial data~\cite{2012gneiting}. Fractal dimension estimation techniques have also been applied for intrinsic dimension analysis in applications where an integer-valued dimension is assumed~\cite{2015camastra,2012mo}.

Several previous studies have observed a relationship between fractal dimension and persistent homology. Estimators based on $0$-dimensional persistent homology (minimum spanning trees) were proposed by Weygart et al.~\cite{1992weygaert} and Martinez et al.~\cite{1993martinez}; we explain their relationship to the current work after Definition~\ref{defn:phdim} below. The PhD thesis of Robins~\cite{2000robins}, arguably the first publication in the field of topological data analysis, studied persistent Betti numbers of fractals and  proved results for the $0$-dimensional homology of totally disconnected sets. In 2012, MacPherson and Schweinhart introduced an alternate definition of $\PH$ dimension, which we refer to here as ``$\PH_i$ complexity'' to avoid confusion (Section~\ref{sec:alternate_phdim}). The two previous notions measure the complexity of a shape rather than a classical dimension; they are trivial for $\mathbb{R}^m,$ for example. A 2019 paper of Adams et al.~\cite{2019adams}, proposed a $\PH$ dimension that is a special case of the one we consider here, and performed computational experiments on self-similar fractals. However, they do not compare the estimates with those of other fractal dimension estimation techniques. In 2018, Schweinhart~\cite{2018schweinhart} proved a relationship between the upper box-counting dimension and the extremal properties of the persistent homology of a metric space, which was a stepping stone to the paper where the current definition was introduced~\cite{2019schweinhart}. The latter paper did not undertake a computational study of the persistent homology dimension; we mention relevant theoretical results in Section \ref{sec:DefOfPHdim}. 

The second class of fractals we study arise in dynamical systems. A dynamical system describes the motion of trajectories $x(t)$ within a phase space, with a rule for describing how $x(t)$ changes as time $t$ increases, and are referred to as either discrete or continuous depending on whether $t \in \Z_+$ or $ t \in \R_+$. 
An attractor $A$ is an invariant set (that is, if $x(t)\in A,$ then $x(s)\in A$ for all $s>t$) in the phase space such that if $x(t)$ is close to $A$, then its distance to $A$ will decrease to zero asymptotically as $ t \to \infty$. Some attractors are simple, consisting of a collection of points or smooth periodic orbits, but some are complicated, such as those shown in Figure \ref{fig:fractals} (d-f). 
In dissipative systems \emph{strange attractors} --- attractors with a fractional dimension --- can occur. 

 For self-similar fractals, dimension is easy to compute exactly and approximate computationally. The situation for strange attractors is quite different; in general, the fractal dimension of an attractor cannot be computed analytically. The measures on these sets may not be regular, and the different notions of dimension appear to disagree in many cases. Note that when estimates for different notions of fractal dimension disagree, it is difficult to tell whether this is because (1) one method converges much faster than the other (2) there is a significant systematic error in the estimates, or (3) the two definitions are genuinely different.
 
 To offer a further cautionary tale,  in \cite{jorba2008mechanism} the authors study a family of quasi-periodically forced $1D$ maps with an attracting invariant curve. As a parameter changes, the curve becomes less stable and it becomes increasingly wrinkled; for any resolution, parameters exist for which the curve will appear to be a strange non-chaotic attractor \cite{grebogi1984strange}. However, until the invariant curve completely loses stability, the attractor remains a smooth 1-dimensional curve. Hence, no matter how fine a resolution we use to compute our dimension estimate, we can never be sure that our estimate is close to its convergent value.

The Lyapunov exponents of a dynamical system describe its exponential rates of expansion and contraction, and can be used to define the Lyapunov dimension of the attractor \cite{young2013mathematical}. 
The Kaplan-Yorke conjecture claims that the  Lyapunov dimension of a generic dynamical system will equal the information dimension of its attractor~\cite{kaplan1979chaotic,frederickson1983liapunov}.  For a class of strange attractors arising from 2D maps, Young showed this to be the case, and moreover that the Lyapunov, Hausdorff, box-counting, and R\'enyi dimensions all coincide  \cite{young1982dimension}. Other attractors appear to exhibit multifractal properties in the sense that different notions of fractal dimension disagree. The multifractal properties of such an attractor can be studied with a one-parameter family of dimensions such as the generalized Hausdorff~\cite{1984grassberger} or R\'{e}nyi dimensions\cite{1987badii,1970renyi}. We can also examine this with the $\PH_i$ dimensions by varying the weight parameter $\alpha$ --- see Figure~\ref{fig:Attractor_spectrum}.

The Takens embedding theorem also illustrates the importance of fractal dimension to the study of attractors. In an experimental setting, it is often impossible to record all of the relevant dynamic variables. 
 However, one can reconstruct  the entire attractor from  the time series of a single observed quantity  using a time delay embedding \cite{1980takens}. 
 Namely, for a choice of time delay $ \tau $ and embedding dimension $ m$,   a 1-dimensional time series 
 $\{ x_1 , \dots x_N \} \subseteq \R $ may be used to construct a $m$-dimensional time series 
  consisting of points $ \{ x_{i}, x_{i-\tau}, \dots , x_{i-(m-1)\tau} \} \in \R^m$. 
 If $m$ is at least twice the box-counting dimension of the attractor, then generically this reconstruction will be diffeomorphic to the original attractor~\cite{sauer1991embedology}. 
 These techniques have been widely applied and we briefly mention some references which use  this in the context of topological data analysis~\cite{garland2016exploring,mischaikow1999construction,myers2019persistent}.

\section{Dimension Estimation Methods}
\label{sec:methods}

\subsection{Persistent Homology Dimension}

\begin{figure}
\centering
\includegraphics[width=.8\textwidth]{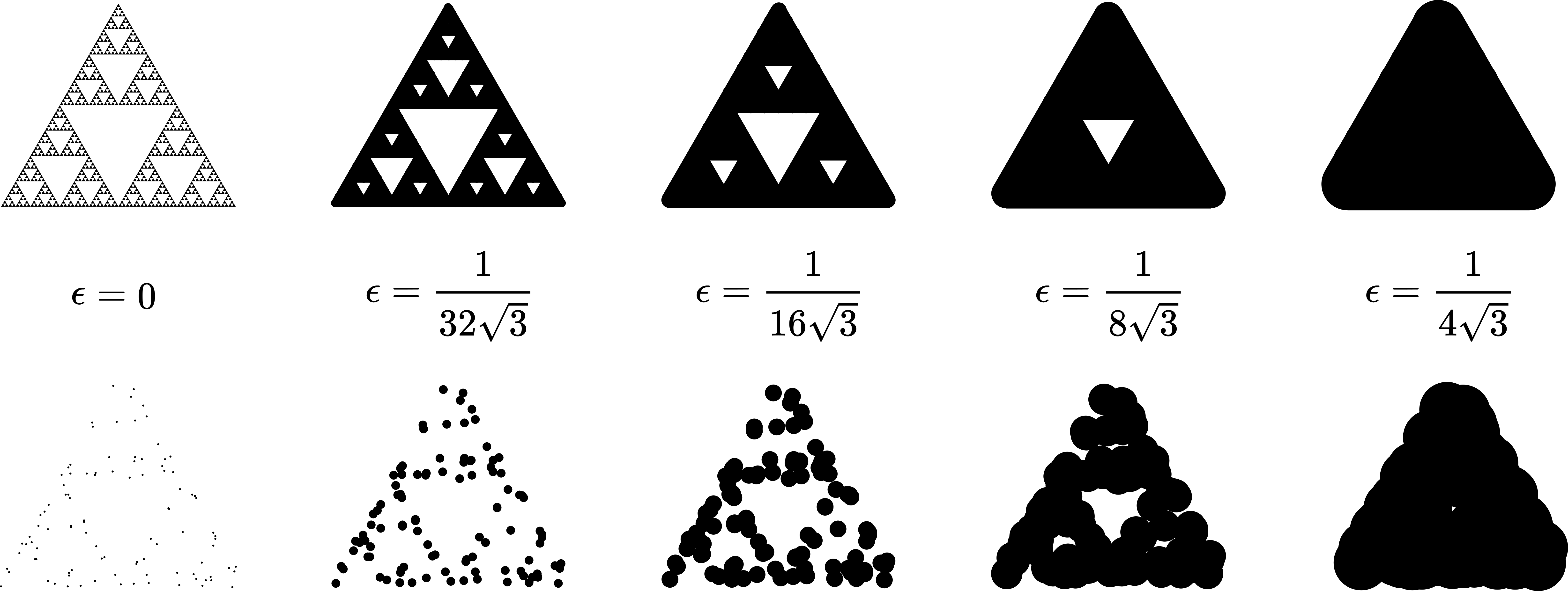}
\caption{\label{fig:epsilonNeighborhoods}
 $\epsilon$-neighborhoods for the Sierpinski triangle (above) and a sample of 100 points from that set (below).}
\end{figure}

\begin{figure}
\centering  
\subfigure[]{
\includegraphics[width=0.3\linewidth]{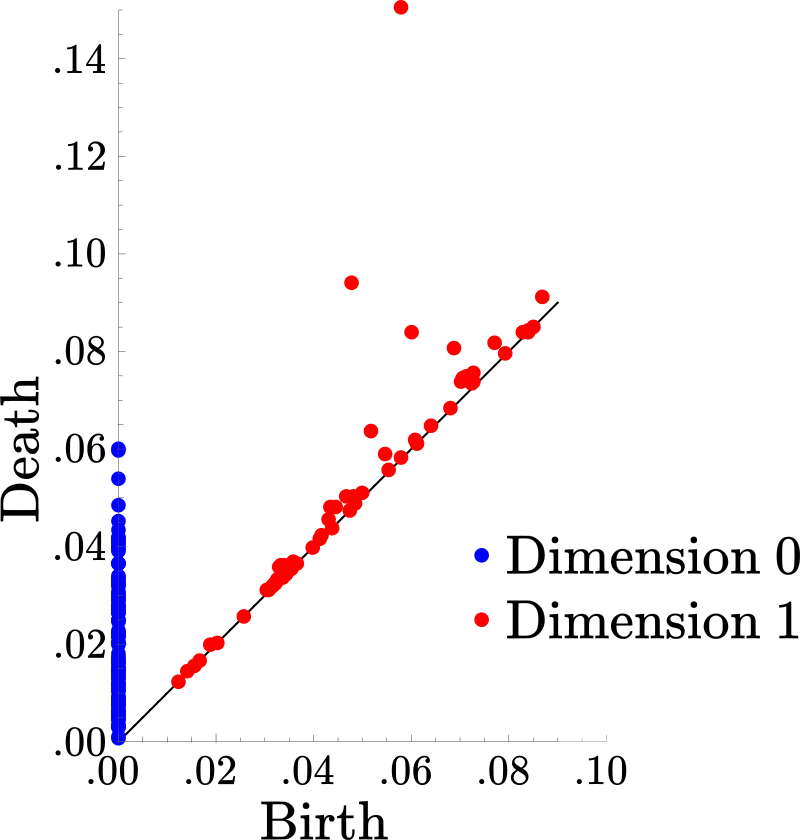}}
\subfigure[]{
\includegraphics[width=0.3\linewidth]{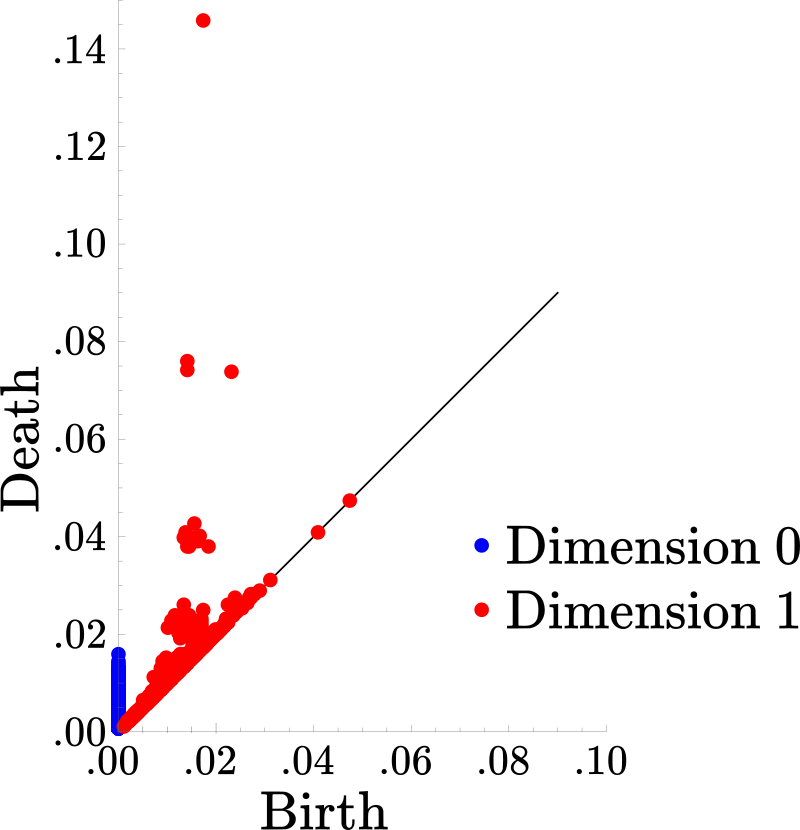}}
\subfigure[]{\label{fig:sierpinski_diagram}
	\includegraphics[width=0.3\linewidth]{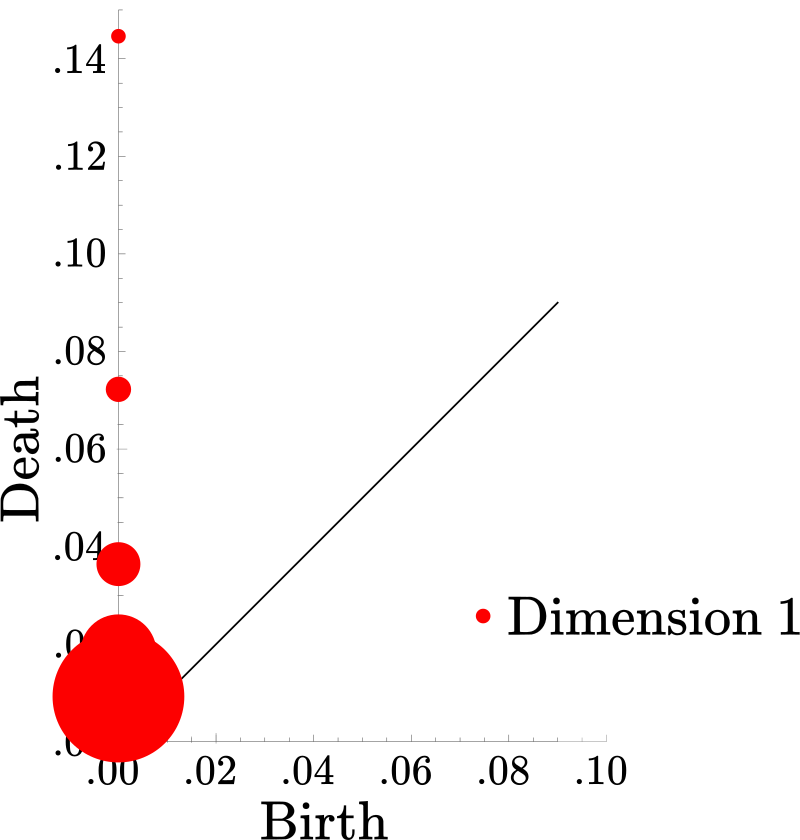}
}
\caption{\label{fig:persistenceDiagrams} Persistence diagrams for (a) 100 point sample from  the Sierpinski triangle, (b) 1000 points from the Sierpinski triangle and (c) the Sierpinski triangle itself (where the area of a dot is proportional to the number of persistent homology intervals with corresponding birth and death time).}
\end{figure}

We give a brief, informal introduction to persistent homology for the special case of a subset of Euclidean space. For a more in depth survey, see~\cite{2008ghrist, 2008edelsbrunner, 2013edelsbrunner, 2016chazal}. Note that the definition of the $\PH$ dimension also makes sense for subsets of an arbitrary metric space.

Persistent homology~\cite{2002edelsbrunner} quantifies the shape of a geometric object in terms of how the topology of the set changes as it is thickened. To be precise, if $X$ is a subset of $\mathbb{R}^m,$ define the family of $\epsilon$-neighborhoods $\set{X_{\epsilon}}_{\epsilon\in\mathbb{R}^+}$ by
\[X_{\epsilon}=\set{y\in \mathbb{R}^m: d(y,X)<\epsilon}\,.\]
Figure~\ref{fig:epsilonNeighborhoods} shows $\epsilon$-neighborhoods of the Sierpinski triangle $S$ and a sample of 100 points from that set. The Sierpinski triangle $S$ contains infinitely many holes which disappear as we thicken it. The first homology group of $S_\epsilon$ is an algebraic invariance which essentially counts the number of holes in $S_\epsilon.$ The last hole to disappear is in the center of the triangle; it vanishes when we have thickened the triangle by $\epsilon=\frac{1}{4\sqrt{3}}.$ Persistent homology represents this largest hole by the single interval $\paren{0,\frac{1}{4\sqrt{3}}}.$ The next three largest holes disappear at  $\epsilon=\frac{1}{8\sqrt{3}}$ and correspond to three intervals  $\paren{0,\frac{1}{8\sqrt{3}}}.$ The following nine holes are represented by nine intervals   $\paren{0,\frac{1}{16\sqrt{3}}},$ and so on.  These are the first dimensional persistent homology intervals of the Sierpinski triangle, which we denote by $\PH_1\paren{S}.$ 

The zero-dimensional persistent homology of a metric space tracks the connected components that merge together as the geometric object is thickened. The Sierpinski triangle is already connected at $\epsilon=0,$ so $\PH_0\paren{S}$ is trivial. However, the finite point sample $\textbf{x}$ shown in Figure~\ref{fig:epsilonNeighborhoods} has 100 components. 
As we thicken $\textbf{x}$ by an amount $\epsilon$ the first component disappears when $\epsilon$ equals $\delta/2,$ where $\delta$ is the smallest pairwise distance between the points. This corresponds to the interval $\paren{0,\delta/2}$ in $\PH_0\paren{\textbf{x}}.$ 
One can find all $\PH_0$ intervals of $\textbf{x}$ by increasing $\epsilon$ and forming an interval whenever there are two points $x_1,x_2\in\textbf{x}$ so that $d\paren{x_1,x_2}<\epsilon$ and $x_1$ and $x_2$ are in different components of $\textbf{x}_{\epsilon_0}$ for all $\epsilon_0<\epsilon.$ 
This is essentially the same as Kruskal's algorithm for computing the minimum spanning tree on $\textbf{x},$ which leads to a proof that there is a bijection between the edges of that tree and the intervals of $\PH_0\paren{\textbf{x}}$ where an interval corresponds to an edge of twice its length.

Roughly speaking, the higher dimensional homology groups $H_i(X_{\epsilon} )$  count the number of higher dimensional ``holes'' in $X_\epsilon.$ The higher dimensional persistent homology is defined in terms of how these groups $H_i\paren{X_{\epsilon}}$ change as $\epsilon$ increases. The structure of the persistent homology is captured by a unique set of intervals that track the birth and death of homology generators as $\epsilon$ changes~\cite{2011cagliaria,2016chazal}. We denote this set of intervals by $\PH_i\paren{X}.$ 

 Traditionally, the information contained in persistent homology is plotted in a persistence diagram showing the scatter plot $\paren{\text{birth},\text{death}}$ for each interval. The persistence diagrams of the Sierpinski triangle, and point samples from it with $100$ and $1000$ points are shown in Figure~\ref{fig:persistenceDiagrams}. Note that the $1$-dimensional persistent homology of the point sample with $1000$ roughly approximates that of the Sierpinski triangle, with $1$ interval that dies around $\epsilon=\frac{1}{4\sqrt{3}},$ three that die around  $\epsilon=\frac{1}{8\sqrt{3}},$ and $9$ that die around  $\epsilon=\frac{1}{16\sqrt{3}}.$ This is a consequence of the ``bottleneck stability'' of persistent homology~\cite{2007cohensteiner}.

It is a remarkable fact that if $\textbf{x}$ is a finite subset of Euclidean space, then the intervals $\PH_i\paren{\textbf{x}}$ can be computed exactly and efficiently. This is done by replacing the infinite family of $\epsilon$-neighborhoods $\set{\textbf{x}_{\epsilon}}_{\epsilon\in\mathbb{R}^+}$ with a finite sequence of finite simplicial complexes, called the Alpha complex of $\textbf{x}$; these complexes are subcomplexes of the Delaunay triangulation on $\textbf{x} $~\cite{2002edelsbrunner}.

Note that there is more than one way to define the persistent homology of a metric space. Here, we use the persistent homology of the \v{C}ech complex persistent homology of $X.$ If $X$ is a subset of Euclidean space, this is equivalent to the persistent homology of the $\epsilon$-neighborhood filtration of $X$ described above (as well as the persistent homology of the Alpha complex, if $X$ is finite). Another common (in-equivalent) notion is the persistent homology of the Vietoris--Rips Complex on $X$~\cite{1927vietoris,2007de_silva}. We use the \v{C}ech complex because there are efficient algorithms to compute its persistent homology for finite subsets of $\mathbb{R}^2$ and $\mathbb{R}^3,$ as discussed below. A persistent homology dimension defined in terms of the Vietoris---Rips complex also has nice properties~\cite{2019schweinhart}, and may work better in examples without an embedding into a small-dimensional Euclidean space.

\subsubsection{Definition of the $\PH$ dimension} 
\label{sec:DefOfPHdim}
Consider a sample $\set{x_1,\ldots,x_n}$ of independent points from the natural measure on the Sierpinski triangle. As $n$ increases, the persistence diagram of the point sample will converge (in the bottleneck distance) to the persistence diagram of the Sierpinski triangle itself. As such, it might seem strange that we can recover the dimension of that set from the $0$-dimensional persistent homology of the samples ($\PH_0(S)$ is trivial). However, one may notice in Figure~\ref{fig:persistenceDiagrams} that there is a large cluster of points, both $0$- and $1$-dimensional, along the plotted diagonal line. 	These points are generally considered to be ``noise'' and do not significantly contribute to the larger features of interest in other applications of persistent homology. However, the rate at which this ``noise'' decays is linked to the dimension of the underlying object. 

To track the growth of this ``noise'', we  define the power-weighted sum for $ \alpha >0$ by: 
\[E_{\alpha}^i\paren{X}=\sum_{I \in \PH_i\paren{X}} \abs{I}^{\alpha}\,,\]
where the sum is taken over all finite intervals and $\abs{I}$ denotes the length of an interval. The scaling properties of random variables of the form $E^0_\alpha\paren{x_1,\ldots,x_n}$ as $n\rightarrow \infty$ has been studied extensively in probabilistic combinatorics~\cite{1992aldous,1996kesten,1988steele,2000yukich}, and the case $i>0$ has recently been of interest~\cite{2018divol, 2019schweinhart}. Motivated by a theorem of Steele~\cite{1988steele}  and the computational work of Adams et al.~\cite{2019adams}, Schweinhart~\cite{2019schweinhart} introduced the following definition of the persistent homology dimension.
\pagebreak[1]
\begin{Definition}
\label{defn:phdim}
Let $X$ be a bounded subset of a metric space and $\mu$ a measure defined on $X$. For each $ i \in \N$ and $ \alpha > 0$ we define the persistent homology dimension:
\[\text{dim}_{\PH_i^\alpha}\paren{\mu}=\frac{\alpha}{1-\beta}\,,\]
where
\[\beta=\limsup_{n\rightarrow\infty} \frac{\log\paren{\mathbb{E}\paren{E_{\alpha}^{i}\paren{x_1,\ldots,x_n}}}}{\log\paren{n}}\,.\]
We write this as the $\PH_i^\alpha$ dimension, and sometimes omit the $i$ or $\alpha$ when making general statements. 
\end{Definition}

That is,  $\text{dim}_{\PH_i^\alpha}\paren{\mu}=d$ if $E_\alpha^i\paren{x_1,\ldots,x_n}$ scales as $n^{\frac{d-\alpha}{d}}.$ Larger values of $\alpha$ give relatively more weight to large intervals than to small ones.   The case $\alpha=1$ is closely related to the dimension studied by Adams et al.~\cite{2019adams}, and agrees with it if defined. Weygaert et al.~\cite{1992weygaert} defined a family of minimum spanning tree dimensions that are equivalent to the $\PH_0$ dimensions, and used heuristic arguments to claim that they coincide with the generalized Hausdorff dimensions for chaotic attractors. Martinez et al.~\cite{1993martinez} asserted that the $\alpha\rightarrow 0$ of the $\PH_0$ dimension gives the Hausdorff dimension for point samples from chaotic attractors.

It is a corollary of Steele~\cite{1988steele} that if $\mu$ is a non-singular measure on $\mathbb{R}^m,$ and $0<\alpha<m$ then $\text{dim}_{\PH_0^\alpha}\paren{\mu}=m.$ Schweinhart~\cite{2019schweinhart} proved that if $\mu$ satisfies a fractal regularity hypothesis called Ahlfors regularity, then $\text{dim}_{\PH_0^\alpha}\paren{\mu}$ equals the Hausdorff dimension of the support of $\mu$ (which coincides with the box-counting dimension under the regularity hypothesis). He also proved that if $d$ equals the upper box-counting dimension of the support of $\mu$ and $\alpha<d$ then  $\text{dim}_{\PH_0^\alpha}\paren{\mu}\leq d,$ as well as weaker results about the cases where $i>0.$

Note that if $\mu$ is supported on a $k$-dimensional subspace of $\mathbb{R}^m,$ then $H_i\paren{X}$ is trivial for $i\geq k.$ It follows that $\text{dim}_{\PH_i^\alpha}\paren{\mu}=0.$ As such, even if $\mu$ is regular, its $\PH_i$ dimension may not equal its Hausdorff dimension unless the Hausdorff dimension is sufficiently large. In particular, if $\mu$ is a $d$-Ahlfors regular measure supported on a $2$-dimensional subspace of $\mathbb{R}^m,$ $d>1.5,$ and $\alpha<d$ then $\PH_1^\alpha$--dimension of $\mu$ equals $d.$~\cite{2019schweinhart}

Schweinhart's results show that the  $\PH_0^\alpha$ and $\PH_1^\alpha$ dimensions of the natural measures on the Sierpinski triangle, Cantor set cross an interval, and Cantor dust equal the Hausdorff dimensions of those sets when $\alpha$ is less than the true dimension.

\begin{figure}
\centering  
\subfigure[]{\label{fig:lorenz_PHfits}
\includegraphics[width=0.3\textwidth]{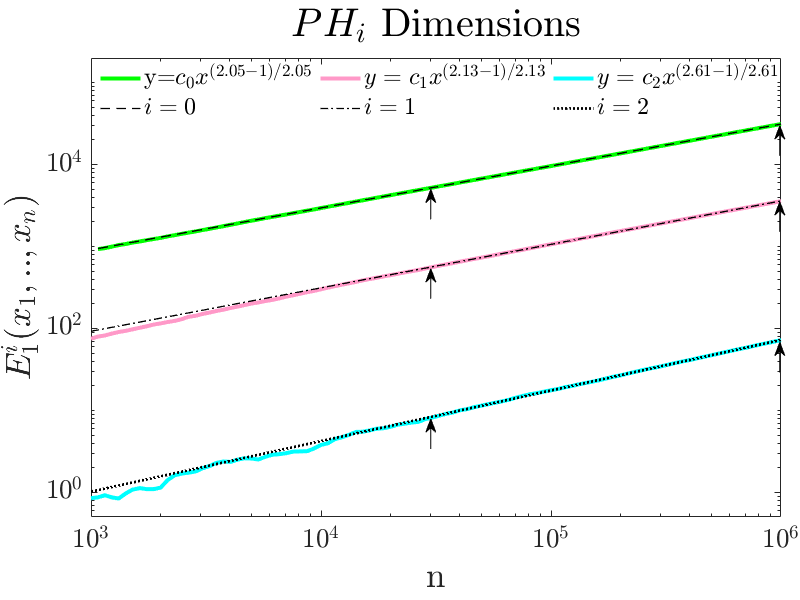}
}
\subfigure[]{\label{fig:LorrenzCorrelationIntegral}
\includegraphics[width=0.3\textwidth]{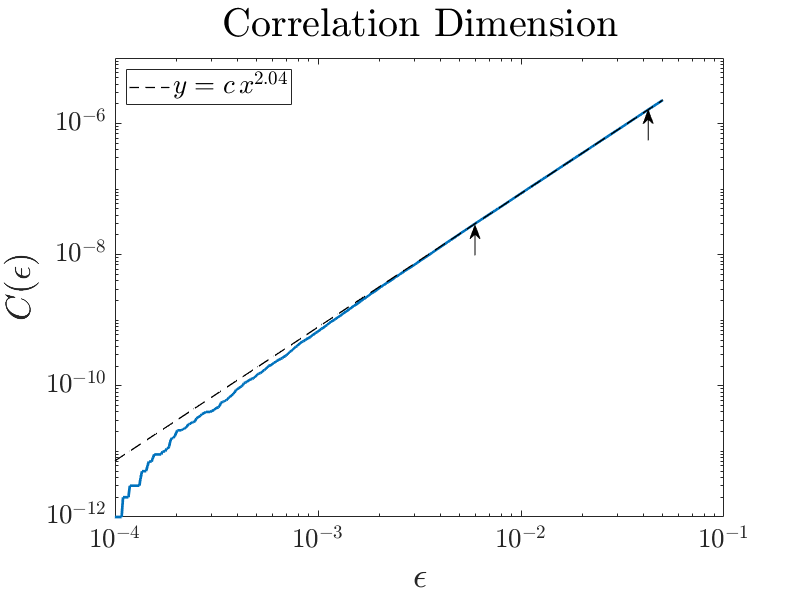}
}
\subfigure[]{\label{fig:LorenzBox}
\includegraphics[width=0.3\textwidth]{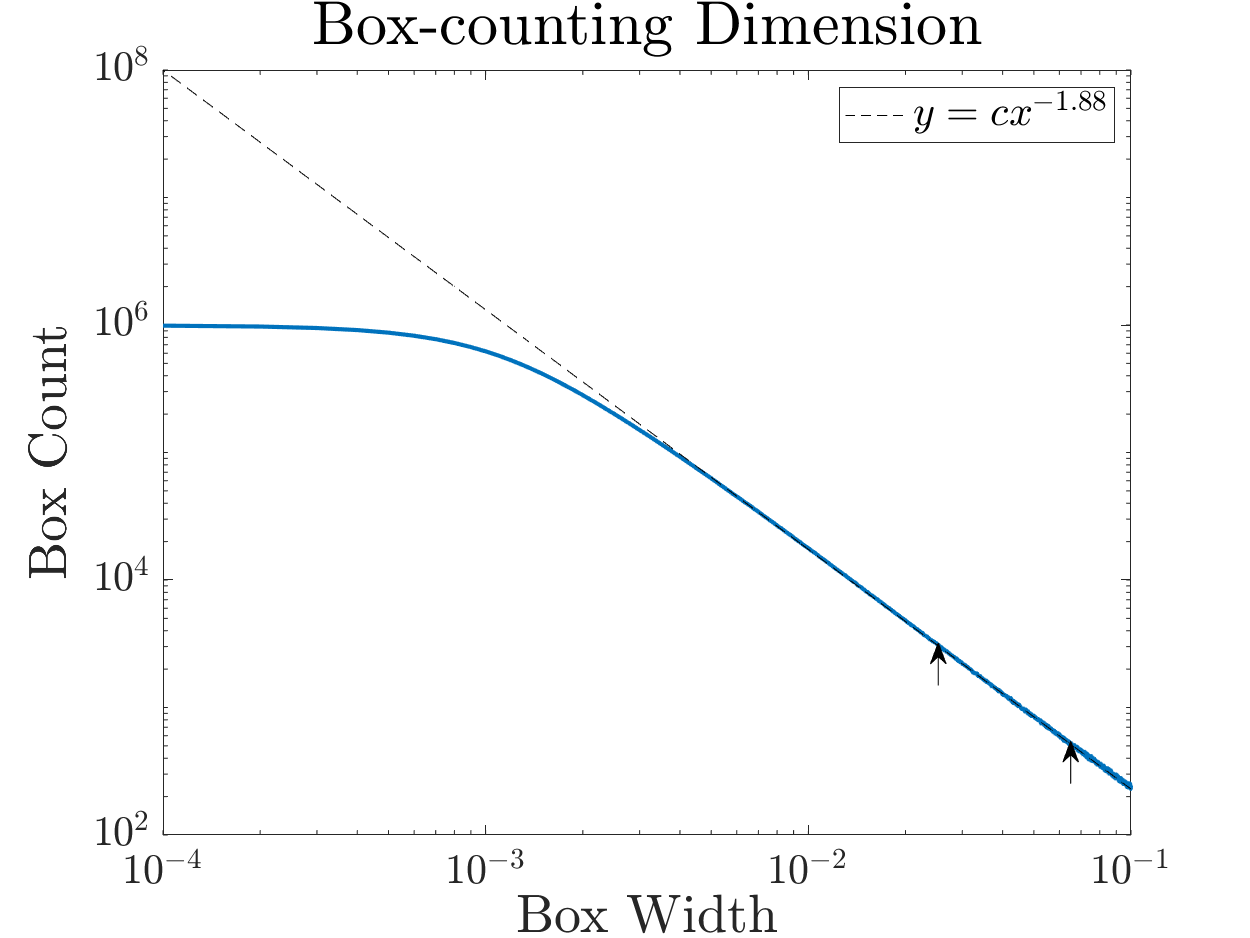}
}
\caption{\label{fig:LorenzFits} Dimension estimation for a sample of $10^6$ points from the Lorenz attractor, using the methods described in Section~\ref{sec:methods}. Power laws were fitted in the ranges denoted by the arrows.}
\end{figure}

\subsubsection{Computation of the $\PH$ dimension}
\label{sec:PH_computation}
We use different methods to compute persistent homology for the cases $i=0$ and $i>0.$ In the $0$-dimensional case, we use minimum spanning tree-based algorithms that work much more quickly, and are fast even for point samples with high ambient dimension. We use the implementation of the dual--tree Boruvka Euclidean minimum spanning tree algorithm~\cite{2010march} included in the mlpack library~\cite{2018curtin} to compute the edges of the minimum spanning tree on a point sample. This algorithm was sufficiently fast for point clouds with $10^6$ points in $\mathbb{R}^8.$ See~\cite{2010march} for a brief survey of other algorithms to compute minimum spanning trees for point sets in Euclidean space and abstract metric spaces.

For $i>0,$ we use GUDHI to compute the persistent homology of the Alpha complex of the point sample.~\cite{GUDHI_alpha,GUDHI_persistence} Both this computation and the previous one could be optimized by re-writing the data structures to support insertion of vertices. 

Given a sample of $n$ points $x_1,\ldots,x_n,$ we compute the $\alpha$-weighted sum $E_\alpha^i\paren{x_1,\ldots,x_{c_i}}$ for $100$ logarithmically spaced values of $c_i$ between $1,000$ and $n.$ Then, we use linear regression to fit a power law to the data $\paren{c_i,E_\alpha^0\paren{x_1,\ldots,x_{c_i}}}.$ After some trial-and-error, we found that fitting the power law between $x=c_n$ and $x=c_{\floor{n}/2}$ provided a reasonable estimate for all examples tested. Using a smaller range sometimes produced better convergence when $i>0,$ but also introduced more oscillations in the estimate. Alternate non-linear regression methods did not seem to result in better performance.

Figure~\ref{fig:lorenz_PHfits} shows the power law fits this method produces for a sample from the Lorenz attractor.

\subsubsection{$\PH$ complexity}
\label{sec:alternate_phdim}
In some instances, the dimension of a metric space can be computed in terms of the persistent homology of the metric space itself. For example, as shown in Figure~\ref{fig:sierpinski_diagram}, the $1$-dimensional persistent homology of the Sierpinski triangle contains intervals that scale as its dimension.   This is captured by an alternate notion of $\PH$ dimension defined by MacPherson and Schweinhart~\cite{2012macpherson}, which measures the complexity of the connectivity of the shape rather than the classical dimension. Here, we refer to it as the ``$\PH$ complexity'' of a shape to differentiate it from notions of dimension. As we will see below, this quantity may be an indicator of when the dimension is hard to estimate using any of the methods presented here. 

If $X$ is a subset of a metric space, define the cumulative $\PH_i$ curve $F_i$ by
\[
F_i\paren{X,\epsilon}
=
\#\left\{ I \in \PH_i\paren{X}:\abs{I}>\epsilon 
\right\}
\]
Then the $\PH_i$ complexity of $X$ is
\[\text{comp}_{\PH_i}\paren{X}=\lim_{\epsilon\rightarrow 0}\frac{-\log\paren{F_i\paren{X,\epsilon}}}{\log\paren{\epsilon}}\,.\]

Note that $\text{comp}_{\PH_i}\paren{\mathbb{R}^n}=0$ for any $i.$ Also, if $S$ is the Sierpinski triangle in Figure~\ref{fig:epsilonNeighborhoods}, our computation of the persistence diagram of $S$ shows that $\text{comp}_{\PH_0}\paren{S}=0$  and $\text{comp}_{\PH_1}\paren{S}=\frac{\log\paren{3}}{\log\paren{2}}.$ 

We can estimate $F_i\paren{X,\epsilon}$   from samples; as the Hausdorff distance between $\set{x_1,\ldots,x_n}$ and $X$ converges to zero, bottleneck stability implies that $F_i\paren{X,\epsilon}$ will converge for values of $\epsilon$ large relative to the Hausdorff distance. In Figure~\ref{fig:ikeda_complexity} below we compute $\text{comp}_{\PH_1}$ of the Ikeda attractor.

\subsection{The Correlation Dimension}

The correlation dimension~\cite{1983grassberger} is commonly used in applications because it is easy to implement and  provides reasonable answers even for relatively small sample sizes.  A probability measure $\mu$ on a metric space $X$ induces a probability measure $\nu$ on the distance set of $X.$ Define the correlation integral of $X$ as the cumulative density function of $\nu$:
\[C\paren{\epsilon}=\mathbb{P}\paren{d\paren{x,y}<\epsilon}.\]

The correlation dimension equals the limit
\[
\lim_{\epsilon\rightarrow 0 } \frac{\log\paren{C\paren{\epsilon}}}{\log\paren{\epsilon}}
\]
if it exists.
There is an extensive literature on the estimation and properties of the correlation dimension; see for example~\cite{1990nerenberg, 1996borovkova,1985takens,1990theiler_b,sprott2001improved}.

\begin{figure} 
\centering  
\subfigure[]{
\label{fig:CD_statTest}
\includegraphics[width=0.45\textwidth]{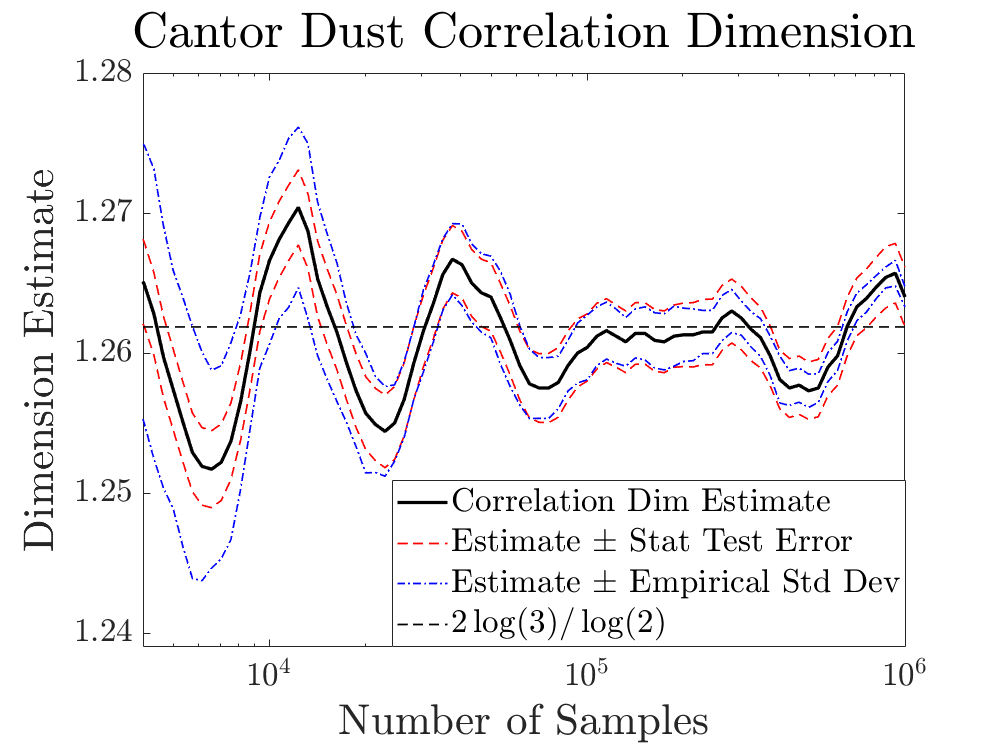}
}
\subfigure[]{
\label{fig:Lorenz_statTest}
\includegraphics[width=0.45\textwidth]{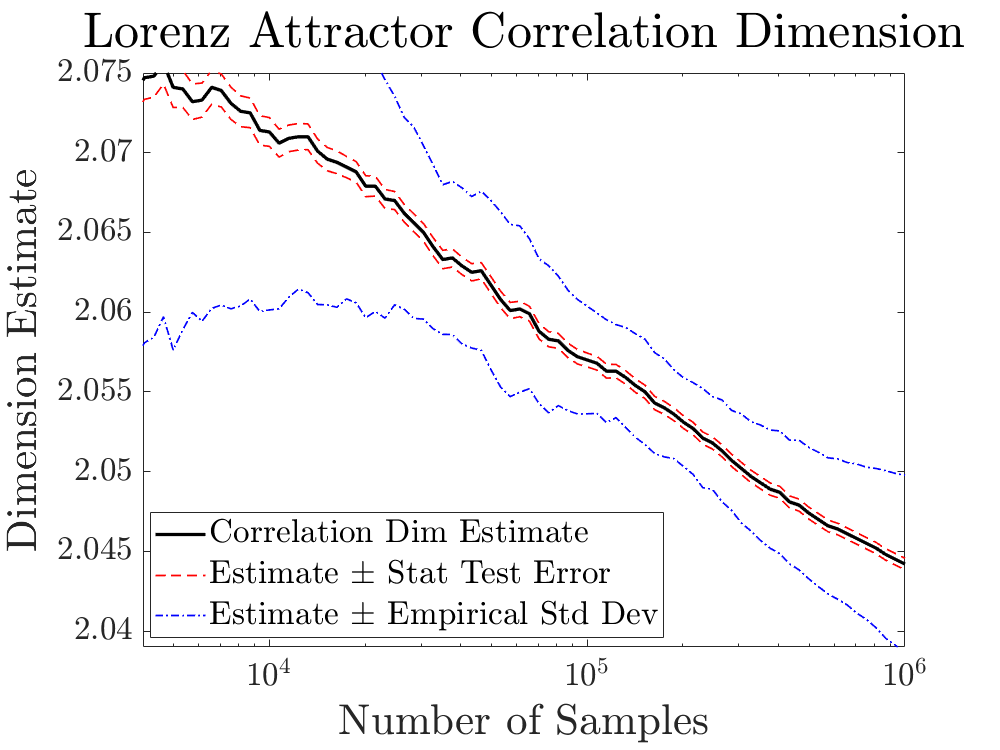}}
		\caption{\label{fig:statistical_test_error} Correlation dimension estimates for (a) the Cantor dust and (b) the Lorenz attractor. The evenly dashed red lines indicate the statistical test error, and the unevenly dashed blue lines show the empirical standard deviation across $50$ trials.}
\end{figure}

\subsection{Computation of the Correlation Dimension}
\label{sec:corrCompute}
For a finite point sample $ x_1, \dots ,x_n$ one can estimate the correlation dimension as 
\[
\lim_{\epsilon \to 0 } \lim_{n \to \infty } \frac{\log C\left(n,\epsilon\right)} {\log \epsilon}
\]
where 
\[
C\left(n,\epsilon\right) = \frac{\# \{ \left(x_i,x_j\right) : d\left(x_i, x_j\right) < \epsilon, i <  j  \}}{n\left(n-1\right)}\,,
\]
if $n$ is taken to $\infty$ appropriately as $\epsilon\rightarrow 0.$  That is, $C\left(n,\epsilon\right)$ measures the number of distances less than $\epsilon$ in proportion to the number of all inter-point distances.

As with other dimension estimates, this limiting expression converges logarithmically slowly. 
To accelerate the convergence, we can essentially apply l'H\^opital's rule to find the limit by computing a slope.  
To do so, we fix a collection of logarithmically spaced values $ \epsilon_1   < \dots < \epsilon_m $, and compute a linear regression through the data $ \left( \log  \epsilon_i ,  \log C(n,\epsilon_i) \right)$.
Our estimate for the correlation dimension is then given by the slope of the line of best fit. 
For a fixed point sample, the values of $ \epsilon_1$ and $\epsilon_m$ are often chosen by hand to avoid outliers and edge effects.  

While computing the $ \mathcal{O}(n^2)$ inter-point distances is prohibative for very large $n$, such a calculation is not needed, as many of these distances are large and do not factor into the dimension calculation \cite{1990theiler}. 
Using a $kd$-tree, one can quickly compute the  $\mathcal{O}(n)$ shortest distances with $\mathcal{O}(n \log n)$ effort.  
After some trial-and-error, we settled upon the  heuristic for choosing $\epsilon_1$ and $\epsilon_m$ as below:
\begin{align*}
 C(N,\epsilon_1) &\approx \frac{n^{.75} }{n(n-1)}&
C(N,\epsilon_m) &\approx \frac{50 n }{n(n-1)}\,.
\end{align*}
See Figure~\ref{fig:LorrenzCorrelationIntegral}. This heuristic provides a way to choose values of $ \epsilon_1$ and $ \epsilon_m$ consistently amongst different fractals and different sample sizes. It 
appears to give near-optimal convergence rates to the true dimension for self-similar fractals, and estimates that agree with the values in the literature for the H\'enon, Ikeda, and Lorenz attractors~\cite{1983grassberger,sprott2001improved,1990theiler}. 

As has been previously reported \cite{1990theiler}, the ``statistical test error'' from the linear regression calculation does not have much predictive value about the limiting dimension --- see Figure~\ref{fig:statistical_test_error}. Also, it is much smaller than the empirical standard deviation of the dimension estimate between trials in some cases.\looseness=-1

\subsection{The Box-counting Dimension}

The box-counting dimension~\cite{1928bouligand}  of a bounded subset of $X$ of $\mathbb{R}^m$ is defined in terms of the number of cubes of width $\delta$ needed to cover $X.$ Let $\set{C_{i}^\delta}_{i\in\mathbb{N}}$ be the cubes in the standard tiling of $\mathbb{R}^m$ by cubes of width $\delta,$ and let $N_{\delta}\paren{X}$ be the number of cubes in $\set{C_{i}^\delta}_{i\in\mathbb{N}}$ that intersect $X.$ Define the upper and lower box-counting dimensions by
\[\text{dim}_{\overline{\text{box}}}\paren{X}=\limsup_{\delta\rightarrow 0}-\frac{N_{\delta}\paren{X}}{\delta}\qquad\text{and}\qquad \text{dim}_{\underline{\text{box}}}\paren{X}=\liminf_{\delta\rightarrow 0}-\frac{N_{\delta}\paren{X}}{\delta}\,,\]
respectively. If the upper and lower box-counting dimensions coincide, the shared value is called the box-counting dimension of $X$ and is denoted  $\text{dim}_{\text{box}}\paren{X}.$ There are several equivalent definitions; see Falconer~\cite{2014falconer} for details. Many studies have investigated the properties and estimation of the box-counting dimension, including~\cite{1994sarkar,1989leibovitch,1991taylor}.

\subsubsection{Computation of the Box-counting Dimension}
\label{sec:boxCompute}

\begin{figure}
\centering  
\subfigure[]{
\includegraphics[width=0.45\textwidth]{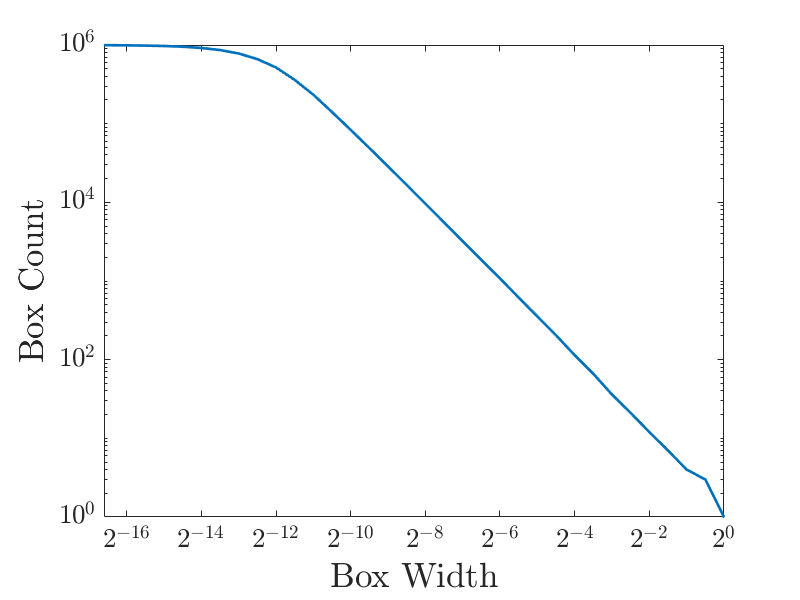}
}
\subfigure[]{
\includegraphics[width=0.417375\textwidth]{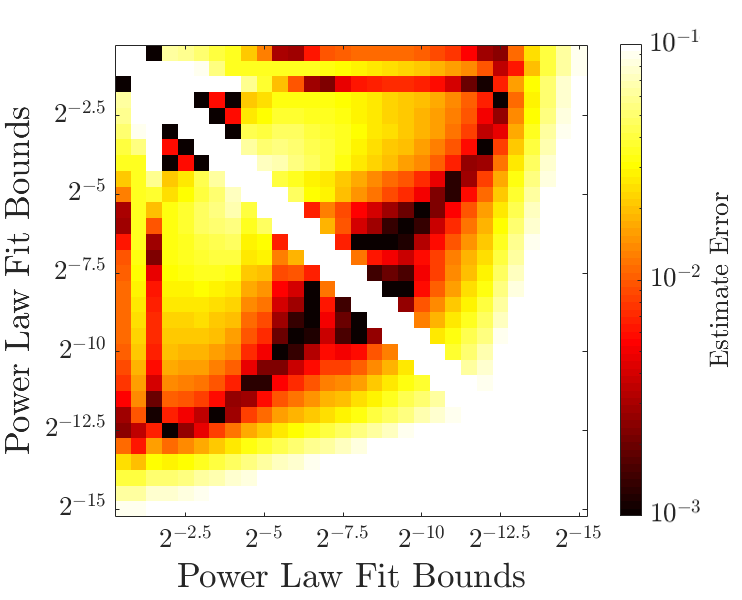}
}
\caption{\label{fig:sierpinski_box} (a) Box-counts for a sample of $10^6$ points from the Sierpinski triangle. (b) Error in the dimension estimate obtained by fitting a power law to the box counts between widths of the form $2^{-i/2}$ and $2^{-j/2}.$ Compare to the error of correlation dimension and $\PH_0^1$ dimension estimates in Figure~\ref{fig:sierpinski2}, which are smaller than .002 for all samples of sizes between $2\,10^5 $ and $10^6$ points. Only  dark squares would yield comparably good estimates for the box-counting dimension.}
\end{figure}

We found it difficult to find a general method to compute the box-counting dimension that worked well for different examples and different numbers of samples. For example, computing box-counts of the form $N_{2^{-i}}\paren{\textbf{x}}$ and fitting a power law between $2^{-i}$ and $2^{-j}$ resulted in estimates that were sensitive to $i,$ $j,$ and the number of samples. In some cases, it was easy to cherry-pick a specific choice just because it seemed to have the best convergence to the true dimension, though other choices resulted in power law fits that looked just as good. See Figure~\ref{fig:sierpinski_box} for an example, which plots the error in the dimension estimate for many possible power law fits. Only very specific choices have errors as small as estimates of the correlation and $\PH_0$ dimensions for the same sample.

We settled on the following method, which produced reasonably good results for planar examples. 
It is based on the observation that if $\set{x_j}_{j\in\mathbb{N}}$ are samples from $X\subset\mathbb{R}^m$  and $\delta>0$ then the box count $N_{\delta}\paren{x_1,\ldots,x_m}$ should converge to $N_\delta\paren{X}$ as $n\rightarrow\infty.$ For a sample $\set{x_1,\ldots,x_m}$, we estimate the box-counting dimension from the smallest boxes for which  $N_{\delta}\paren{x_1,\ldots,x_m}$ appears to have stabilized. 

 For a a family of point samples $\textbf{x}_m$ in $\mathbb{R}^n$ (where the sizes of $\textbf{x}_m$ are logarithmically spaced as before) we rescale and translate the point samples so they are contained in a unit cube. Then, we compute box-counts of the form $N_{i/10000}\paren{\textbf{x}_m}$ for $1\leq i \leq 1000.$ We fit a power law to the data $\paren{i/10000,N_{i/10000}\paren{\textbf{x}_m}}$ in the range $(\left \lceil{j/2}\right \rceil 
,j),$ where $j$ is the smallest index so that 
\[N_{j/10000}\paren{\textbf{x}_{m}}=\text{min}_{1\leq k \leq 4}\paren{N_{j/10000}\paren{\textbf{x}_{m-k}}}\,.\] We used linear regression to fit the power law; non-linear regression did not produce substantially different results. See Figure~\ref{fig:LorenzBox}.

We tried several variations. For example fitting the power law in the range $(\left \lceil{.9 j}\right \rceil ,j)$ produced estimates that converged faster with $n$ for some examples but exhibited large oscillations in others.

\section{Results for Self-similar Fractals}
\label{sec:examples_discrete}
We compare the performance of the fractal dimension estimation procedures for four different self-similar fractals: the Sierpinski triangle ($S$), the Cantor dust ($C\times C$), the Cantor set cross an interval ($C\times I$), and the Menger sponge ($M$). The first three are subsets of $\mathbb{R}^2$, and $M$ is contained in $\mathbb{R}^3.$ We chose these sets to illustrate the observed relationship between $\text{comp}_{\PH_0}\paren{X}$ and the performance of dimension estimation techniques; see Table~\ref{table:frac}. Definitions of the sets and sampling methods are contained in Appendix~\ref{appendix:self_similar}.

\begin{table}
\begin{tabular}{l|c|c|c|c|c|c}
Example & True Dim. &  $\text{comp}_{\PH_0}\paren{X}$ & $\text{comp}_{\PH_1}\paren{X}$ &  $\text{comp}_{\PH_2}\paren{X}$ \\
\hline
$S$ & $\frac{\log(3)}{\log(2)} $  & $0$  &  $\frac{\log(3)}{\log(2)} $& $-$\\
\hline
$C\times I$  & $1+\frac{\log(2)}{\log(3)} $  &  $\frac{\log(2)}{\log(3)} $  & $0$ & $-$ \\
\hline
$C\times C$  &  $\frac{2\log(2)}{\log(3)} $  &  $\frac{2\log(2)}{\log(3)} $  &  $\frac{2\log(2)}{\log(3)}$ &  $-$  \\
\hline
$M$ &  $\frac{\log(20)}{\log(3)} $  &  $0$ &  $\frac{\log(20)}{\log(3)}$ & $\frac{\log(20)}{\log(3)}$
\end{tabular}
\caption{\label{table:frac} Data for the four self-similar fractals studied here. As discussed in Section~\ref{sec:alpha}, lower values of $\alpha$ appear to yield better dimension estimates when $\text{comp}_{\PH_i}\paren{X}\neq 0.$}
\end{table}
 
{
\begin{table}
\begin{tabular}{l|c|c|c|c|c|c|c|c|c}
 & True & Correl. & Box & $\PH_0^{.5}$ & $\PH_0^{1}$ & $\PH_1^{.5}$ & $\PH_1^1$ &  $\PH_2^{1}$ & $\PH_2^2$  \\
\hline
$S$ & $\approx 1.585$  & $1.585 $ & $1.586 $ & $1.585  $  &  $1.585  $ & $1.587  $ &   $1.620  $ & $-$ & $-$  \\
\hline
$C\times I$  & $\approx 1.631$  &  $1.633  $ & $1.618  $ & $1.629  $  &  $1.623  $ & $1.634 $ & $1.642  $  & $-$ & $-$ \\
\hline
 $C\times C$ &  $\approx 1.262$   &  $1.263  $ & $1.267  $ & $1.263  $  &  $1.289  $ & $1.268  $ & $1.303  $& $-$ & $-$    \\
\hline
$M$ & $\approx 2.727$  & $2.716$  & $2.703$ & $2.705$ &  $2.706$  & 2.878 & $2.773$ & $2.945$ & $2.881$
\end{tabular}
\caption{\label{table:self_similar_estimates} Dimension estimates for self-similar fractals, averaged over $10$ trials of $10^6$ samples.}
\end{table}
}

\begin{figure}
\centering  
\subfigure[]{
\label{fig:sierpinski1}
\includegraphics[width=0.46\textwidth]{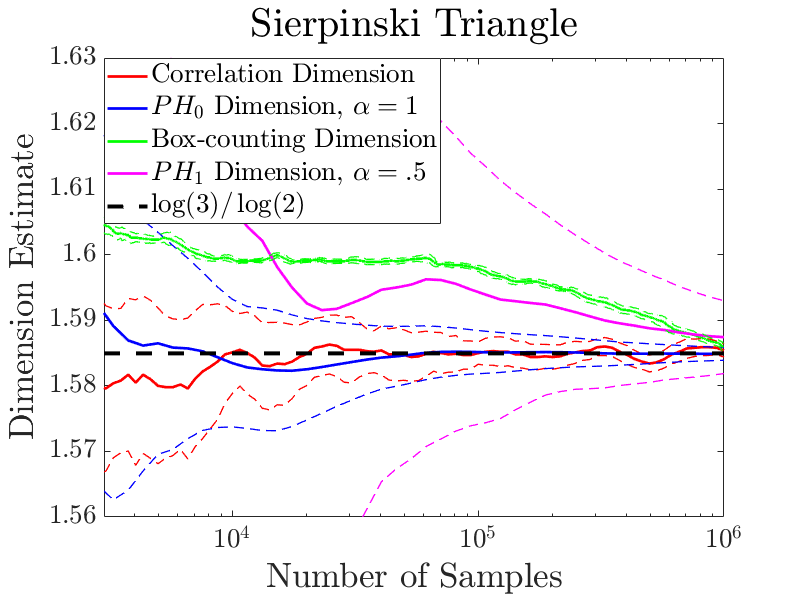}
}
\subfigure[]{
\label{fig:sierpinski2}
\includegraphics[width=0.46\textwidth]{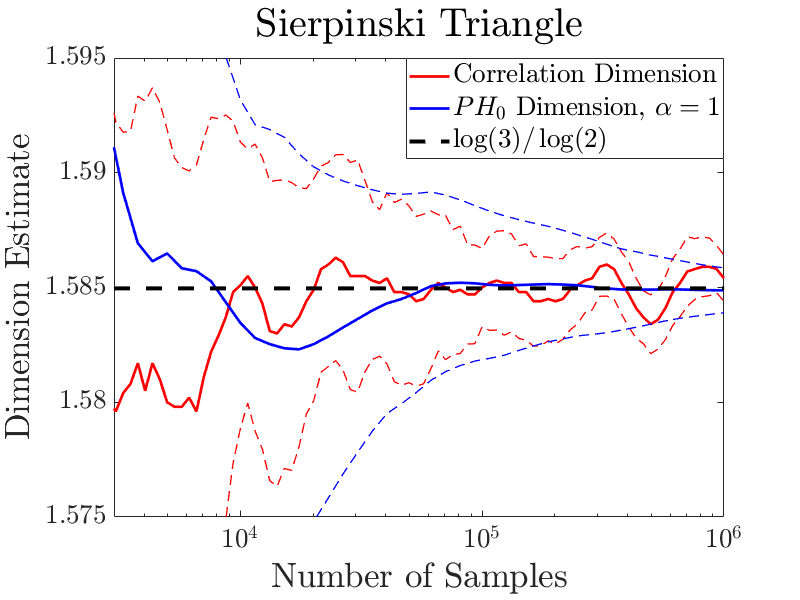}}
\subfigure[]{
\label{fig:cxi_comparison}
\includegraphics[width=0.46\textwidth]{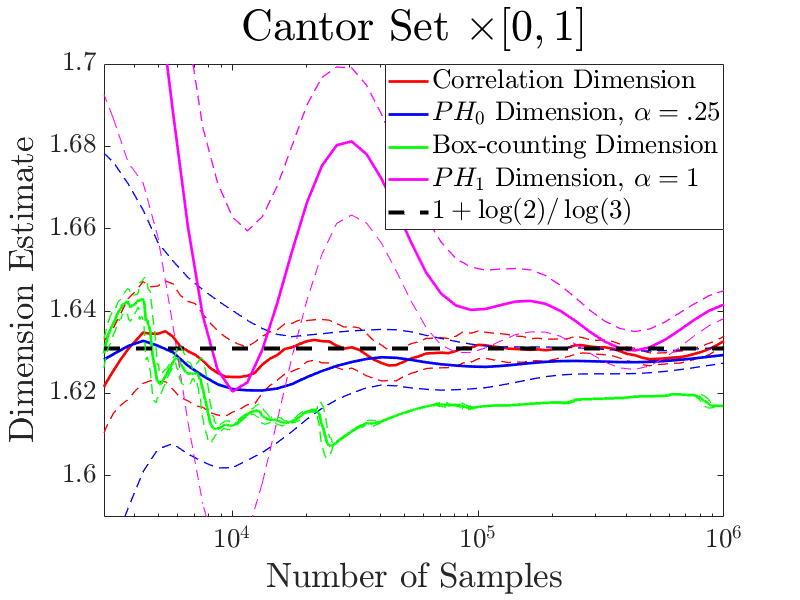}
}
\subfigure[]{
\label{fig:cd_comparison}
\includegraphics[width=0.46\textwidth]{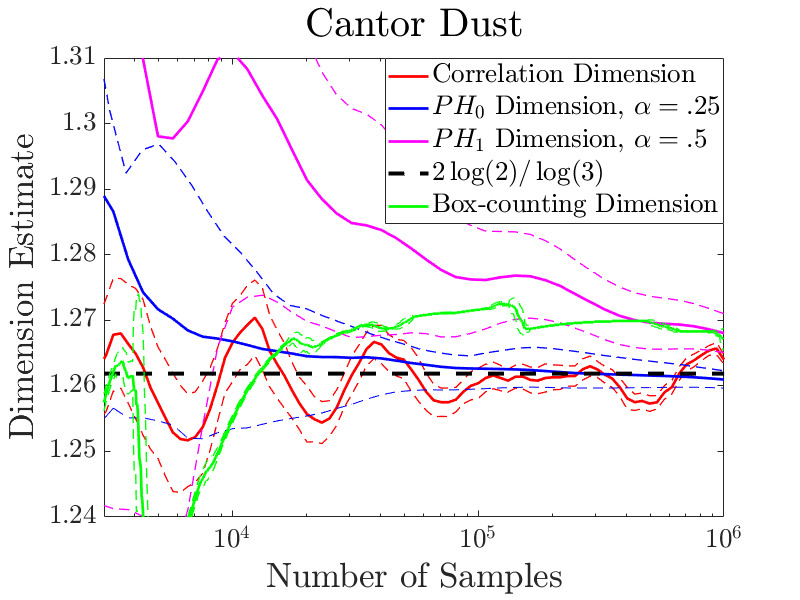}}
\subfigure[]{
\includegraphics[width=0.46\textwidth]{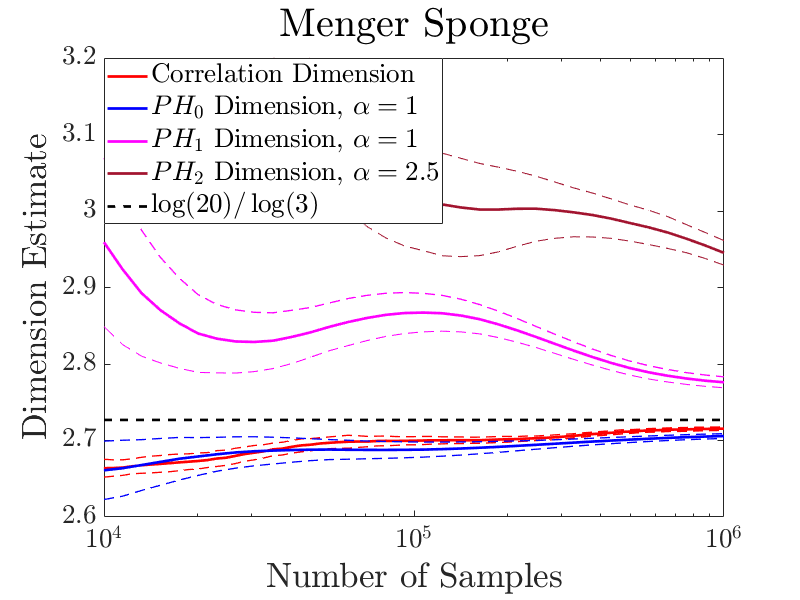}
}
\caption{\label{fig:self_similar_estimates} Dimension estimates for self-similar fractals. Note that (b) is a close-up of (a). We omitted the box-counting dimension in (e) because the automated procedure described in Section~\ref{sec:boxCompute} did not result in stable estimates until $n$ was larger than $10^5.$}
\end{figure}

For each set,  we sample $50$ trials of $10^6$ points from the corresponding natural measure. We compute 100 dimension estimates for each trial, at logarithmically spaced numbers of points between $10^3$ and $10^6.$  To show the convergence of the dimension estimate to the true value, we plot number of points against the dimension estimate averaged across the trials, with thinner dotted lines one standard deviation above and below the estimate.

Dimension estimates for the four examples are plotted in Figure~\ref{fig:self_similar_estimates}. In all examples, the $\PH_0$ dimension and the correlation dimension perform better than the $\PH_1$ or box-counting dimensions. The relatively poor performance of the $\PH_1$ dimension can be ascribed to the fact that the weighted sums $E_{\alpha}^1\paren{x_1,\ldots,x_n}$ are smaller and noisier than the corresponding sums $E_{\alpha}^0\paren{x_1,\ldots,x_n}$ (see Figure~\ref{fig:lorenz_PHfits}). As mentioned previously, it was difficult to find an effective general method to produce box-counting estimates. This is illustrated for the Sierpinski triangle in Figure~\ref{fig:sierpinski_box}. While the plot of box width vs. box count looks linear for a large range on the log--log plot, only very specific choices of bounds will produce dimension estimates with errors on the same order as the $\PH_0$ and correlation dimensions. These bounds vary unpredictably with the sample size and example.

The convergence of the $\PH_0$ and correlation dimensions to the true dimension is very similar in all four examples. 
The correlation  dimension estimates exhibit oscillations as the sample size increases, a phenomenon often ascribed to the fractal's lacunarity~\cite{1990theiler}. The $\PH_0$ dimension estimate appears less sensitive to oscillations, but this is likely due to how the power law fit is performed in the computations (as described in Section~\ref{sec:PH_computation}); choosing a narrower range over which to fit a power law results in a more oscillatory estimate. That said, for sizable samples, it is impractical to compute the correlation integral over a large range. The $\PH_0$ dimension provides a computationally practicable way to use information from multiple lengthscales. It also has the advantage of having the parameter $\alpha$ which can be tuned to give a better estimate --- see Section~\ref{sec:alpha} for a discussion on the choice of $\alpha.$ The variance of the correlation dimension estimate is slightly lower, but this doesn't mean much when oscillations are present.\looseness=-1

\begin{figure}
\centering  
\subfigure[]{
\label{fig:correlationperformance}
\includegraphics[width=0.3\linewidth]{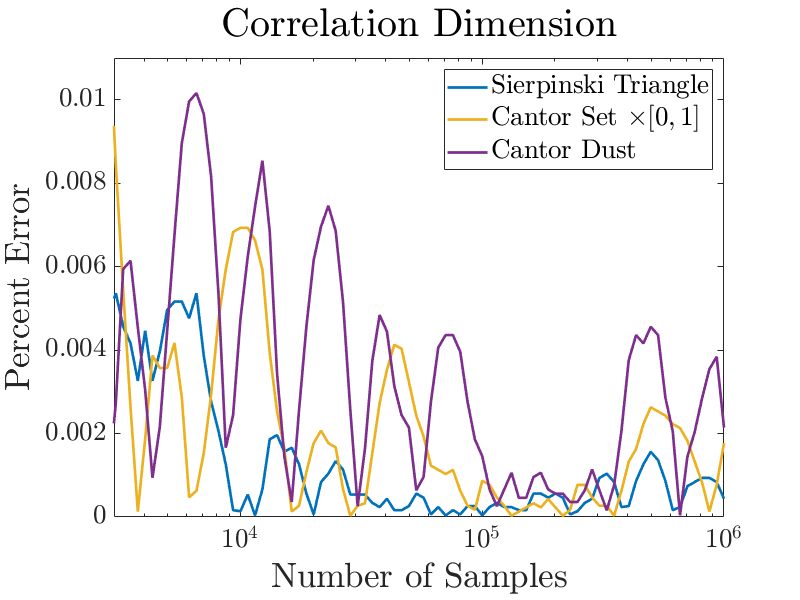}
}
\subfigure[]{
\label{fig:ph0performance}
\includegraphics[width=0.3\linewidth]{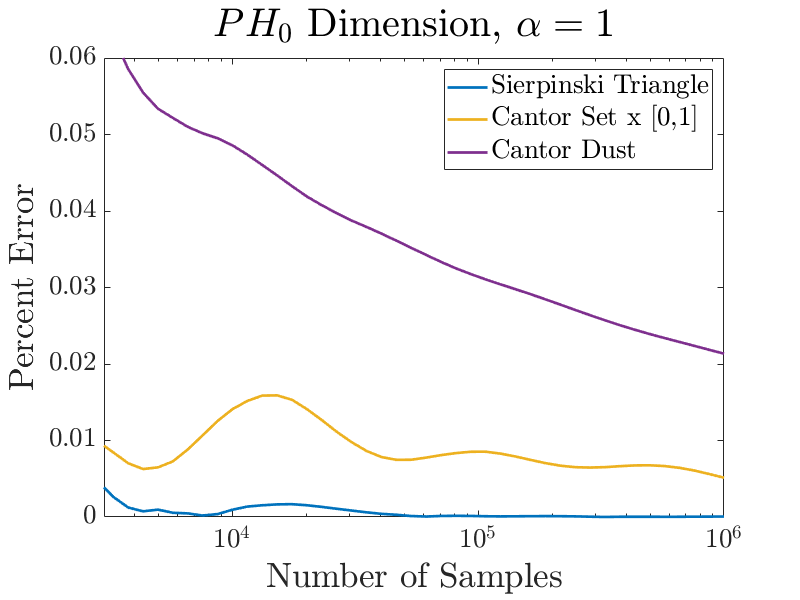}}
\subfigure[]{
\label{fig:ph1performance}
\includegraphics[width=0.3\linewidth]{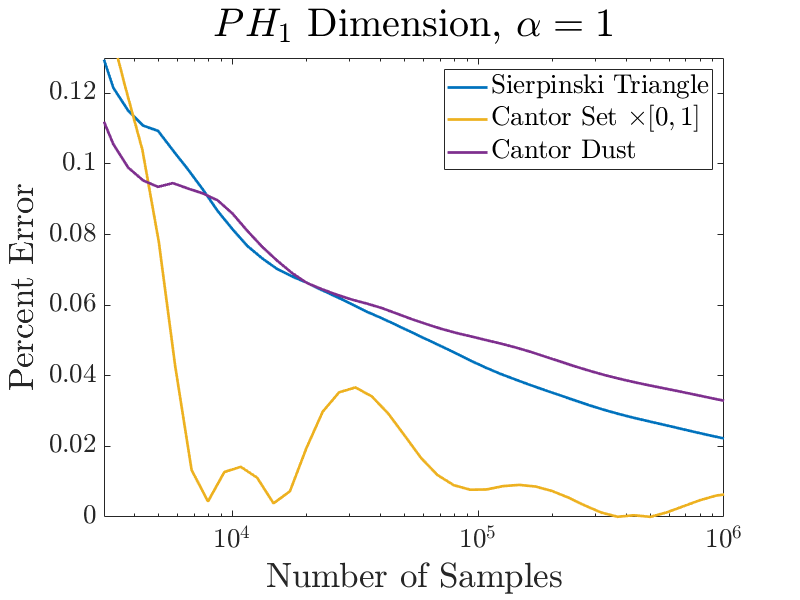}}
		\caption{\label{fig:error}Percentage error of dimension estimates for three self-similar fractals.}
\end{figure}

\subsection{Percentage Error}

An interesting pattern emerges when we compare the percentage error of the dimension estimates across the three planar examples, as in Figure~\ref{fig:error}. Both $\PH_0$ and correlation dimensions perform best for the Sierpinski triangle ($\text{comp}_{\PH_0}(S)=0$), second best for the Cantor set cross an interval ($\text{comp}_{\PH_0}(C\times I)=\text{dim}\paren{C\times I}-1=\frac{\log\paren{2}}{\log\paren{3}}),$ and worst for the Cantor dust ($\text{comp}_{\PH_0}(C\times C)=\text{dim}\paren{C\times C}=\frac{2\log\paren{2}}{\log\paren{3}}$). That is, dimension estimation appears to be more difficult when the connectivity of the underlying set is more complex. Note that $\text{comp}_{\PH_0}(X)$ can be computed (see~\cite{2012macpherson}), and could be used as an indicator of whether additional caution is warranted when discussing dimension estimation results.

For the $\PH_1$ dimension, the situation is different and the difficulty of dimension estimation appears to (unsurprisingly) be related to $\text{comp}_{\PH_1}$ rather than $\text{comp}_{\PH_0}.$ The rate of convergence was fastest for the Cantor set cross an interval ($\text{comp}_{\PH_1}(C\times I)=0),$ but slower for the Cantor dust  ($\text{comp}_{\PH_1}(C\times C)=\text{dim}\paren{C\times C}=\frac{2\log\paren{2}}{\log\paren{3}}$) and Sierpinski triangle  ($\text{comp}_{\PH_1}(S)=\text{dim}\paren{S}=\frac{\log\paren{3}}{\log\paren{2}}$).

We exclude the Menger sponge from these figures, as the ambient dimension likely influences the difficulty of dimension estimation. 

\begin{figure}
\centering  
\subfigure[]{\label{fig:sierpinski_ph0_alpha}
\includegraphics[width=0.375\textwidth]{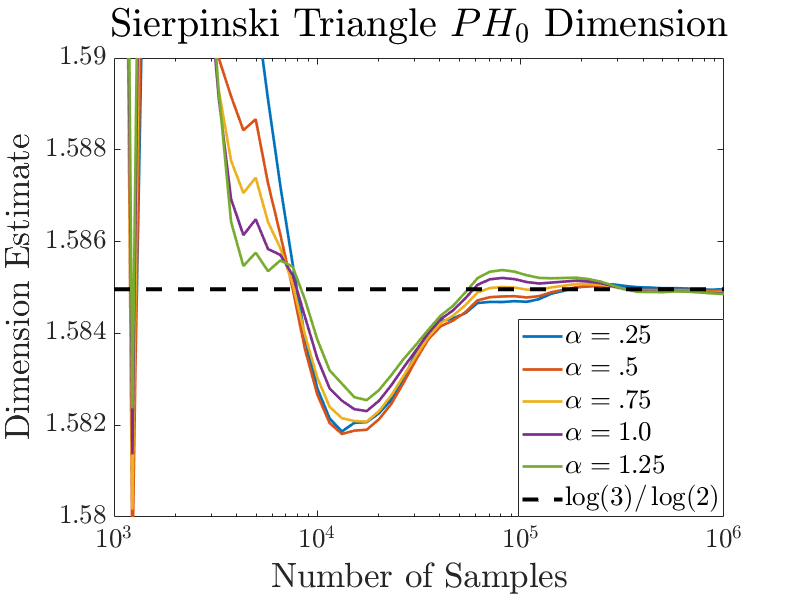}
}
\subfigure[]{\label{fig:sierpinski_ph1_alpha}
\includegraphics[width=0.375\textwidth]{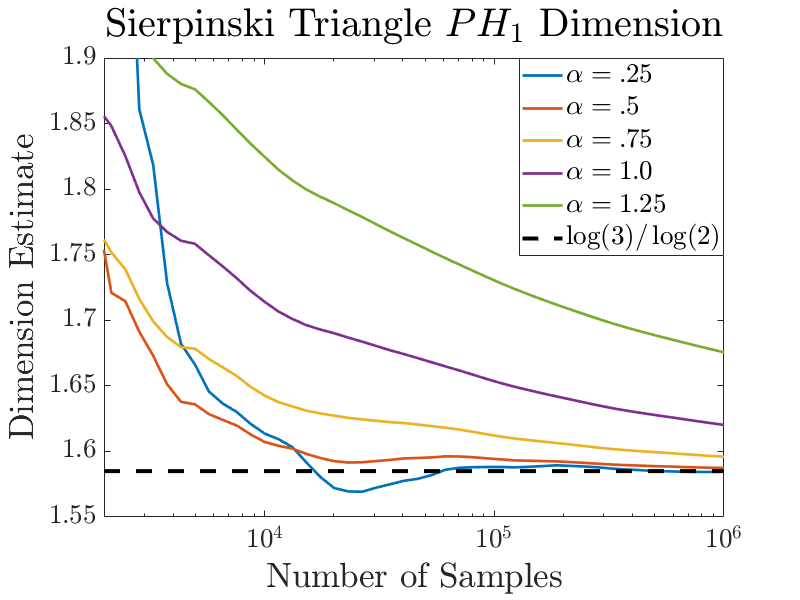}
}
\subfigure[]{\label{fig:cxi_ph0_alpha}
\includegraphics[width=0.375\textwidth]{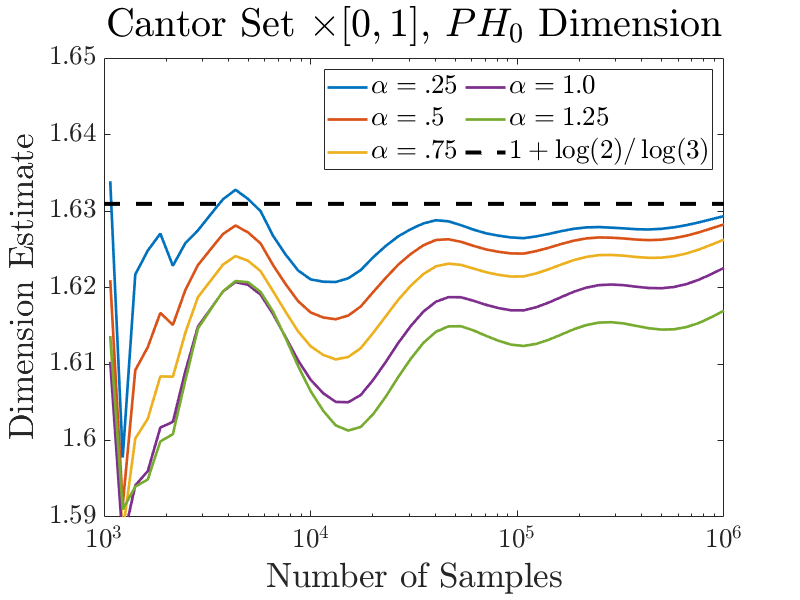}
}
\subfigure[]{\label{fig:cxi_ph1_alpha}
\includegraphics[width=0.375\textwidth]{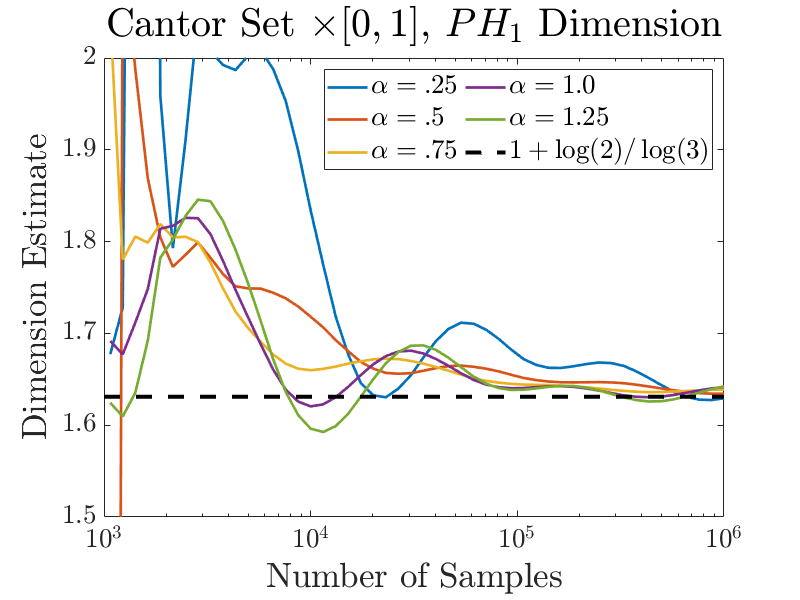}
}
\subfigure[]{\label{fig:cd_ph0_alpha}
\includegraphics[width=0.375\textwidth]{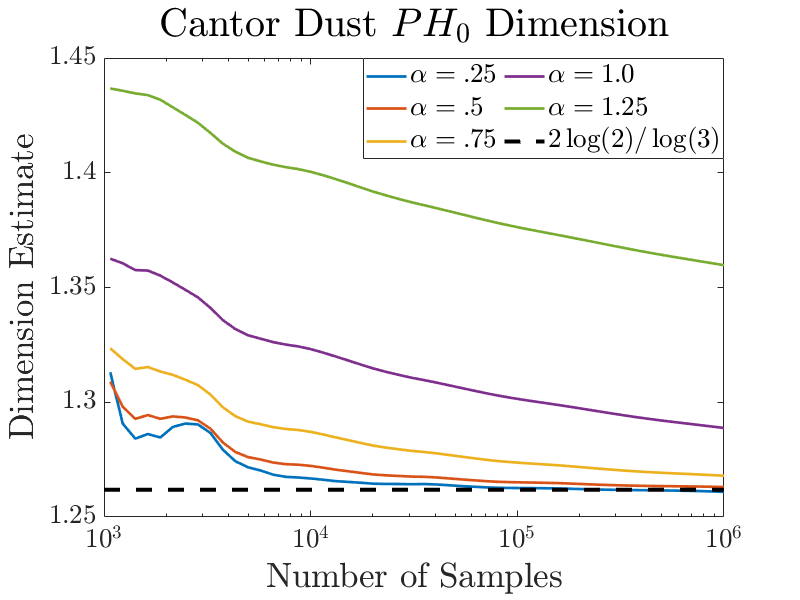}
}
\subfigure[]{\label{fig:cd_ph1_alpha}
\includegraphics[width=0.375\textwidth]{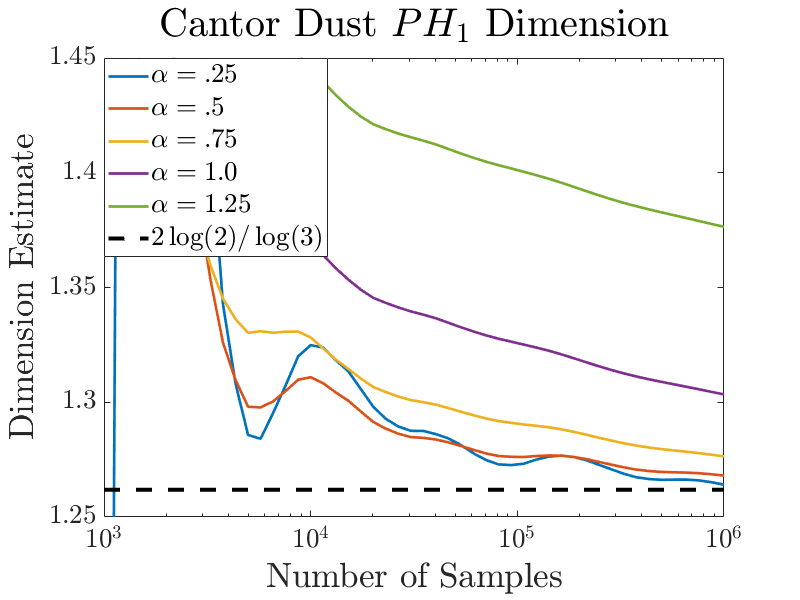}
}
\subfigure[]{\label{fig:sierpinski_alpha_variance_a}
\includegraphics[width=0.375\textwidth]{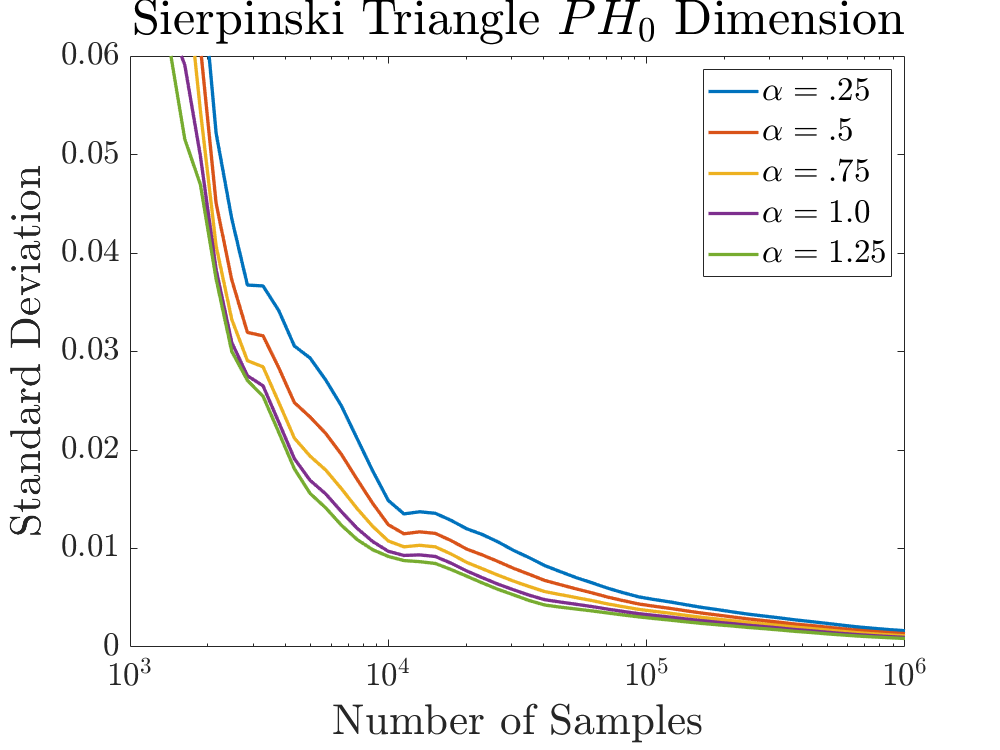}
}
\subfigure[]{\label{fig:sierpinski_alpha_variance_b}
\includegraphics[width=0.375\textwidth]{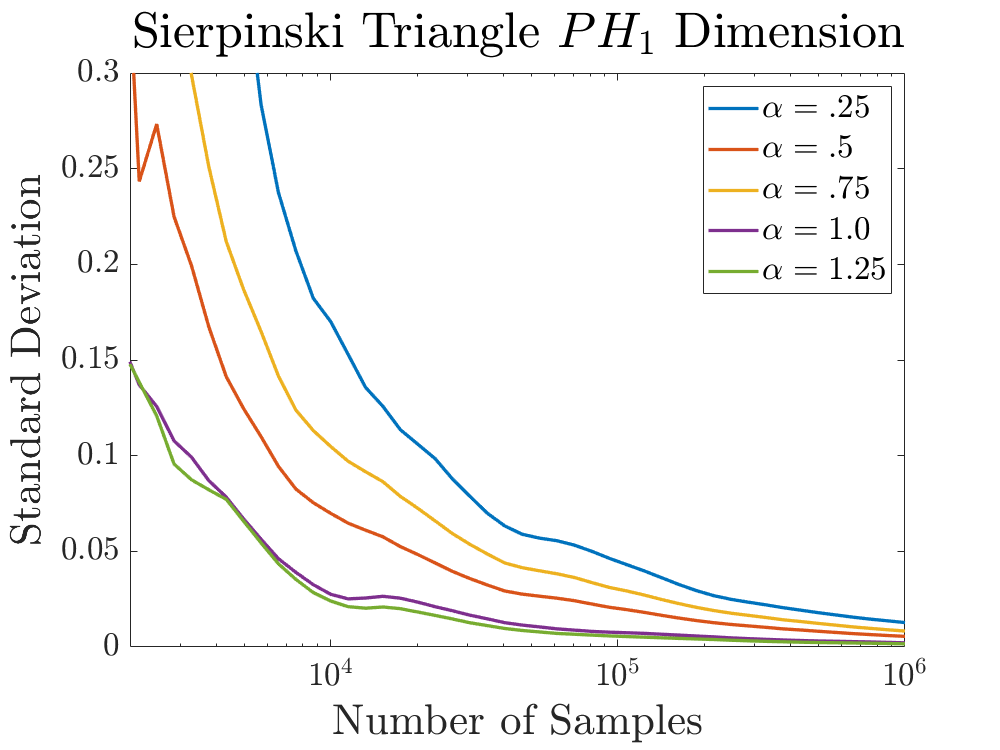}
}
\caption{\label{fig:alpha} $\PH_i$ dimension estimates for various choices of $\alpha.$ (g) and (h) show the standard deviation of the dimension estimate between trials.}
\end{figure}

\subsection{Dependence on $\alpha$}
\label{sec:alpha}

Estimates of the $\PH_0^\alpha$ and $\PH_1^\alpha$ dimensions for various choices of $\alpha$ are shown in Figures~\ref{fig:alpha}. In the cases where $\text{comp}_{\PH_i}\paren{X}=0,$ there are only small differences between dimension estimates for different choices of $\alpha,$ with perhaps a slight advantage for higher values of $\alpha$ (Figures~\ref{fig:sierpinski_ph0_alpha} and~\ref{fig:cxi_ph1_alpha}).  
However, when   $\text{comp}_{\PH_i}\paren{X}>0,$ dimension estimates for different values of $\alpha$ are substantially different. Lower values of $\alpha$ yield better estimates when $\text{comp}_{\PH_i}(X)>0$ for planar examples, and middle values of $\alpha$ (i.e. equal to about half the true dimension) seem to provide the best convergence for $i=1,2$ (but convergence is slow, especially for $i=2.$) Lower values of $\alpha$ appear to give dimension estimates that have a higher variance between samples --- see Figures~\ref{fig:sierpinski_alpha_variance_a} and~\ref{fig:sierpinski_alpha_variance_a} for the Sierpinski triangle, which is representative.

\section{Results for Attractors}
\label{sec:examples_dynamics}

\begin{figure}
\centering  

\includegraphics[width=0.45\textwidth]{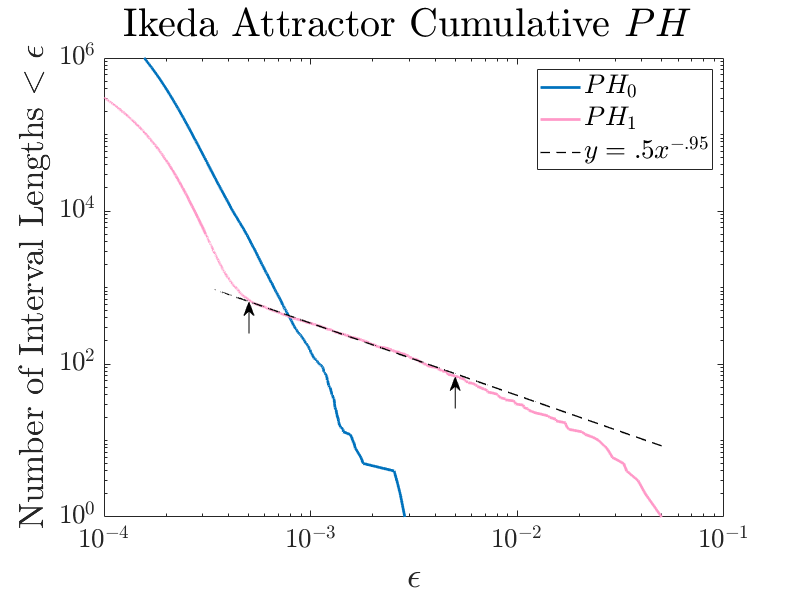}
		\caption{\label{fig:ikeda_complexity} To estimate $\text{comp}_{\PH_1}$ of the Ikeda attractor, we fit a power law in the range marked by the arrows. Note that $\PH_1$ exhibits two regimes: noise from the point sample and persistent features of the attractor itself. $\PH_0$ is only noise, as the Ikeda attractor is connected.}
\end{figure}

\begin{table}[ht]
	\begin{tabular}{l|c|c|c|c|c|c}
		& Correlation & Box-counting  & $\PH_0^{.5}$ & $\PH_0^{1}$ & $\PH_1^{.5}$ & $\PH_1^1$  \\
		\hline
		H\'enon   & $1.21 $ & $1.22  $ & $1.28  $  &  $1.28  $ & $1.52  $ &   $1.36  $  \\
		\hline
		Ikeda   &  $1.68  $ & $1.72  $ & $1.71  $  &  $1.72  $ & $1.71  $ & $1.72  $   \\
		\hline
		Rulkov   &  $1.01  $ & $1.52  $ & $1.62  $  &  $1.87  $ & $2.02$ & $<2.13$   \\
		\hline
		Lorenz   &  $2.04  $ & $>1.90  $ & $2.06 $  &  $2.05  $ & $<2.14$ & $<2.12$   \\
		\hline
		MG   &  $3.04 $ & $>2.45  $  & $3.59  $  &  $3.70  $ &  $-$ &  $-$
	\end{tabular}
	\caption{\label{table:attractor_estimates} Dimension estimates for chaotic attractors, averaged over $10$ trials of $10^6$ samples. }
\end{table}

\begin{figure}
\centering  
\subfigure[]{
\label{fig:henon_comparison}
\includegraphics[width=0.47\textwidth]{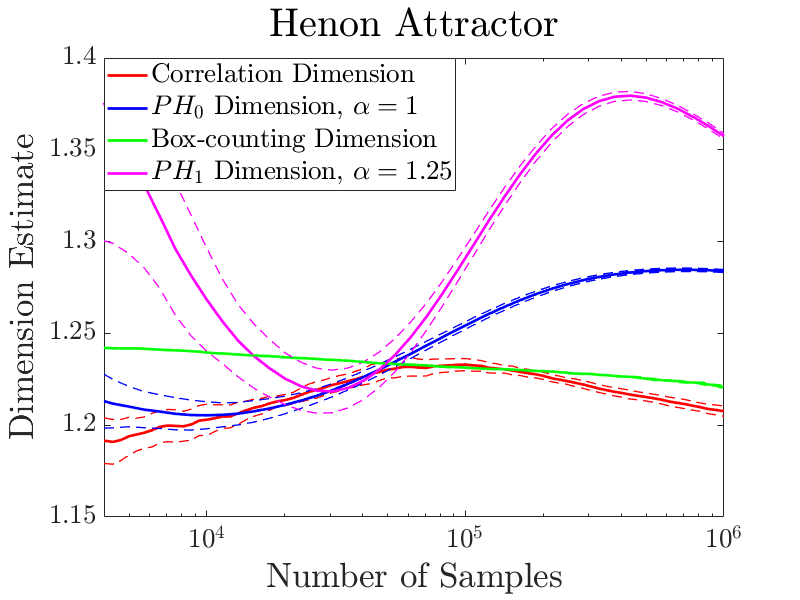}
}
\subfigure[]{
	\label{fig:ikeda_comparison}
	\includegraphics[width=0.47\textwidth]{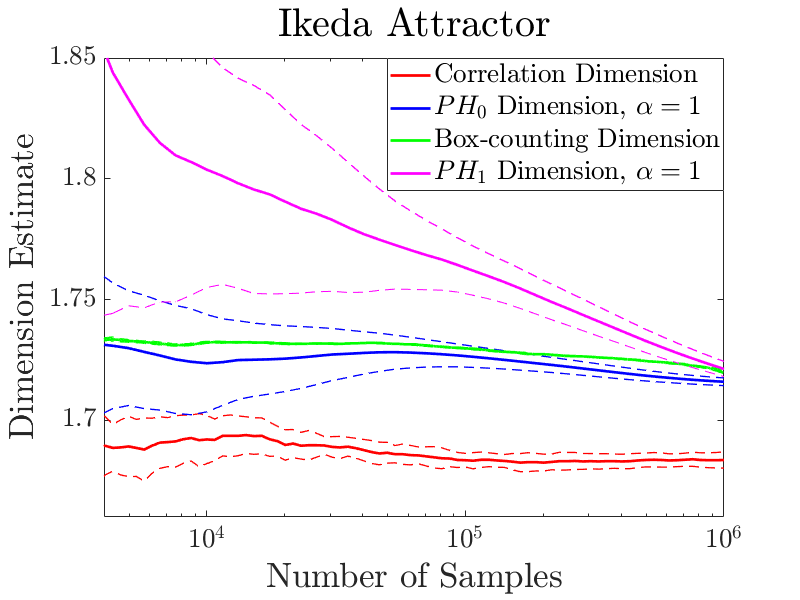}
}
\subfigure[]{
	\label{fig:rulkov_comparison}
	\includegraphics[width=0.47\textwidth]{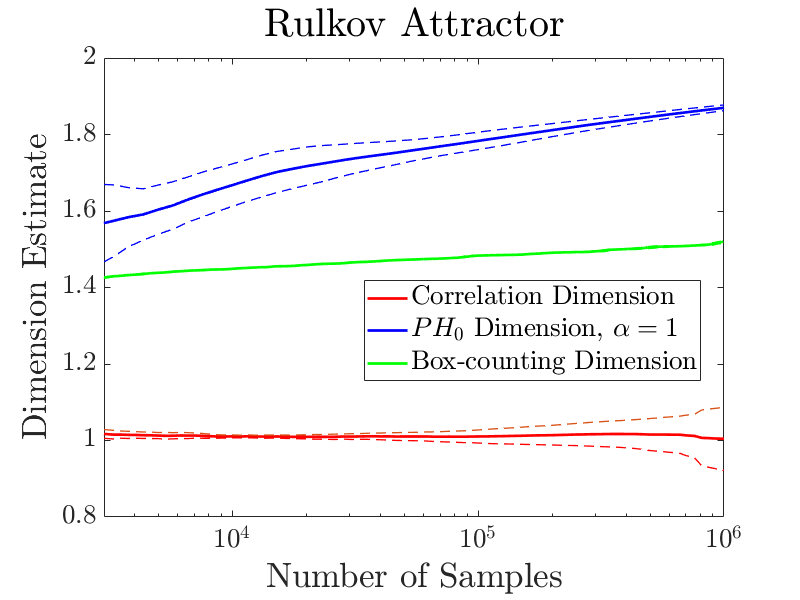}
}
\subfigure[]{
\label{fig:lorenz_comparison}
\includegraphics[width=0.47\textwidth]{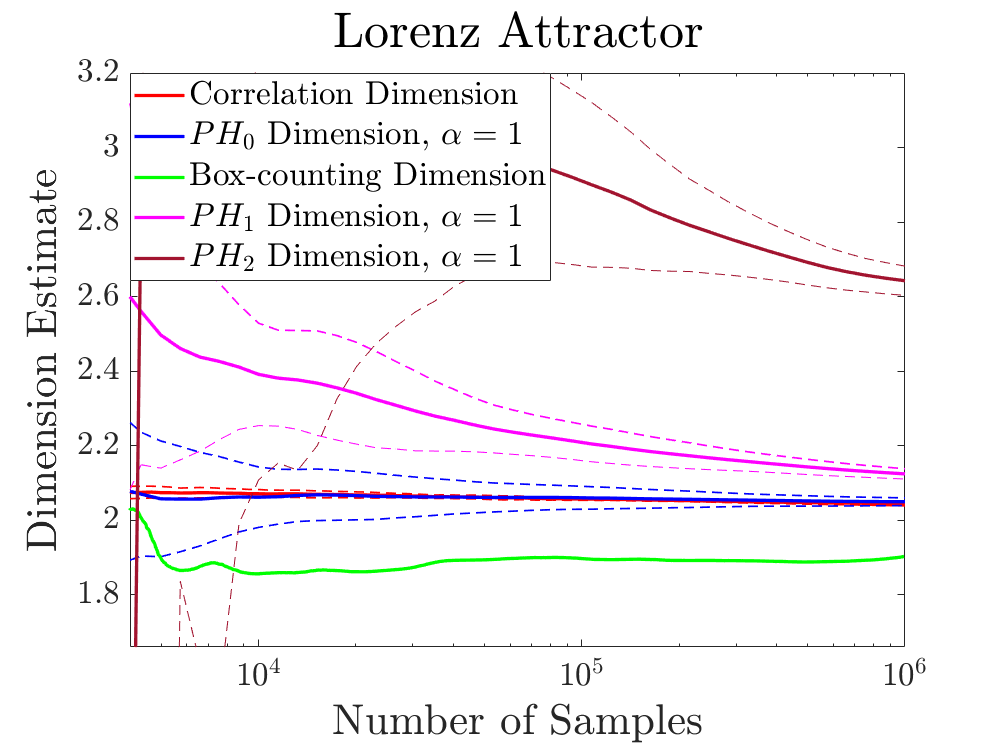}
}
\subfigure[]{
	\label{fig:MG_comparison}
	\includegraphics[width=0.47\textwidth]{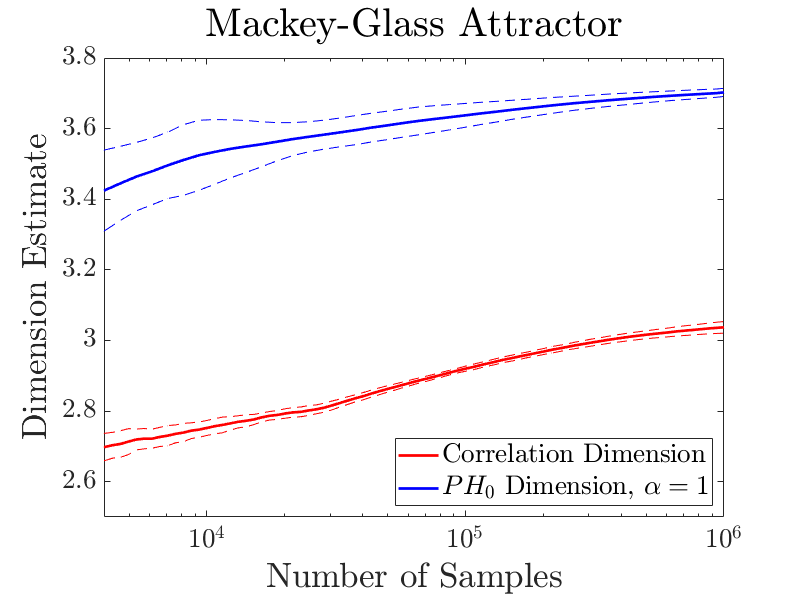}
}

		\caption{\label{fig:Attractor_data}Average dimension estimates for various attractors, with dashed lines denoting $ \pm 1$ sample standard deviation. Note the different scales.}
\end{figure}

For each chaotic attractor we sample $50$ randomly selected initial conditions, and generate a time series  of $10^6$ points after discarding an  initial transient trajectory (except the Mackey-Glass attractor, where we sample $25$ initial conditions). For each time series, we compute 100 dimension estimates for each trial, at logarithmically spaced numbers of points between $10^3$ and $10^6.$ 
The various dynamical systems we studied are described in Appendix \ref{sec:Attractors_Description} and the dimension estimates we obtained for $ N = 10^6$ are summarized in Table \ref{table:attractor_estimates} and Figure \ref{fig:Attractor_data}. 

In all cases, we are skeptical that the box-counting dimension has approached a limiting value. The upper box dimension is known to be an upper bound for the $\PH_0^\alpha$~\cite{2019schweinhart} and correlation dimensions, but most of our estimates of these dimensions tend to be higher than our estimates of the box-counting dimension. This comports with previous observations that box-counting estimates do not seem to converge well for strange attractors~\cite{1982greenside}.

\begin{figure}
	\centering  
	\subfigure[]{
		\label{fig:rulkov_ph0}
		\includegraphics[width=0.45\linewidth]{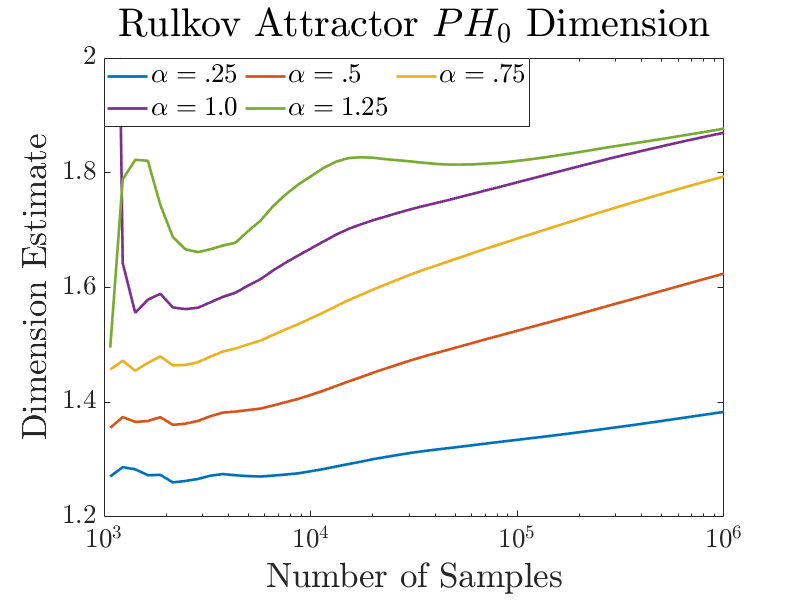}
	}
	\subfigure[]{
		\label{fig:MG_ph0}
		\includegraphics[width=0.45\linewidth]{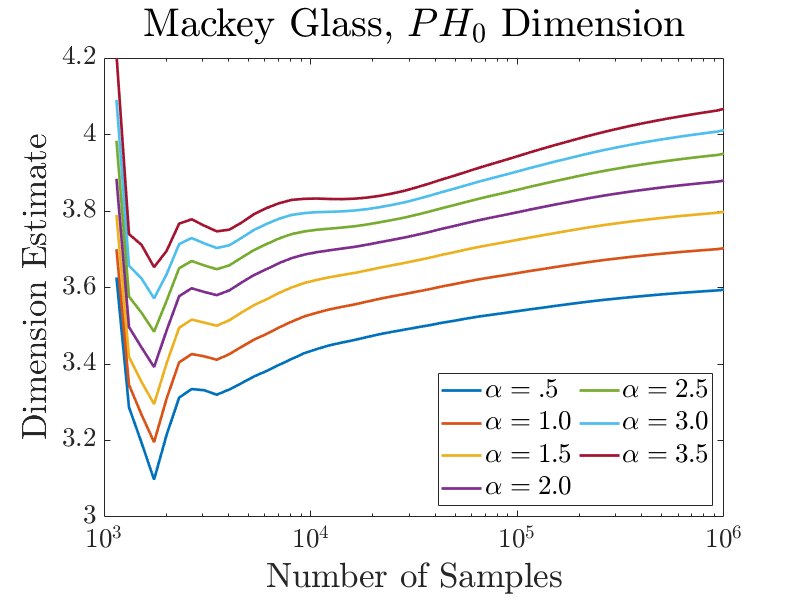}}
	\caption{\label{fig:Attractor_spectrum} Attractor $\PH_i$ dimension estimates for various choices of $\alpha.$ }
\end{figure}

Of the two-dimensional maps we studied, the Ikeda map appears to be most well-behaved.  The box-counting, $\PH_0$, and even the  $\PH_1$ dimension all seem to converge to a common value of $1.71$, with a correlation dimension of $1.68$ not that far off. The Ikeda attractor $I$ was the only chaotic attractor we studied with non-trivial persistent homology at multiple lengthscales (rather than the persistent homology ``noise'' of point samples). Here, we compute its $\PH_1$ complexity, $\text{comp}_{\PH_1}\paren{I}.$ Figure~\ref{fig:ikeda_complexity} shows the cumulative interval count function defined in~\ref{sec:alternate_phdim}, for a sample of $10^6$ points from the Ikeda attractor. The $\PH_0$ curve only shows very small intervals, which come from the noise of point samples. This is as expected, because the Ikeda attractor is connected. However, the $\PH_1$ curve shows two different regimes, one that corresponds to noise and one that appears to be picking up on the elongated holes visible in Figure~\ref{fig:Ikeda}. The cumulative length plot appears to follow a power law in the range $\epsilon=.0005$ to $\epsilon=.005.$ Fitting a power law to this range gives an estimate of  $\text{comp}_{\PH_1}\paren{I}=.95.$ This is lower than the dimension estimates for the Ikeda attractor, indicating that the elongated holes in that set scale at a different rate than its other properties.

For the H\'enon attractor, both $\PH_0$ and $\PH_1$ dimension estimates had large oscillations which make it hard to say whether they will converge to values different from either the box-counting or correlation dimensions. The correlation dimension exhibited smaller oscillations, which comports with previous observations by Theiler \cite{theiler1988quantifying},  who also warns of the danger of estimating the dimension of  attractors with long-period oscillations. Previous studies~\cite{1987arneodo} have claimed that the H\'enon attractor is multifractal for the parameters chosen here, but we are not confident enough in our dimension estimates to make such an assertion.

The Rulkov attractor is non-homogeneous, as easily evidenced by Figure \ref{fig:Rulkov}, and not surprisingly our fractal dimension estimates differed from each other, ranging from $ 1.00$ (correlation dimension) to $2.13$ ($\PH_1^1$ dimension). 
The $\PH_1^\alpha$ estimates began above $2$ for small samples sizes and decreased toward $2$ as the number of samples increased. 
For various values of $\alpha,$ the $\PH_0^\alpha$ dimension estimates ranged from $1.38$ for $\alpha=1$ to $1.88$ for $\alpha=1.25$ (see Figure \ref{fig:Attractor_spectrum}). As one may see in Figure \ref{fig:rulkov_comparison}, the variance of the correlation dimension estimate increases with sample size,  suggesting the heuristic we chose for fitting the correlation integral requires further fine-tuning for this example. 
The Rulkov map is noninvertible, so it does not belong to the  class of 2D maps considered in \cite{young1982dimension} where many of the fractal dimensions are known to coincide. While it is difficult to determine if numerical dimension estimates have converged, our results certainly suggest that the various fractal dimension definitions may not agree for the Rulkov attractor. 

For the Lorenz attractor,  the $\PH_0$ and the correlation dimensions performed well and both appear to converge toward $2.05$. However, neither the $\PH_1$ nor the $\PH_2$ dimensions perform well. In fact, nearly $ 10^5$ are needed before the $\PH_2$ dimension estimate less than 3, the ambient dimension in which the points reside! As in other cases, this is likely due to the fact that higher dimensional persistent homology contains less information --- see Figure~\ref{fig:LorenzFits}. The sum of the lengths of the $\PH_1$ and $\PH_2$ intervals of a point sample of $10^6$ points are $10$ and $400$ times smaller than the sum of the lengths of the $\PH_0$ intervals, respectively. Also, the plot of $n$ versus $E_1^i\paren{x_1,\ldots,x_n}$ is much noisier for $\PH_2$ than $\PH_0.$

Lastly we consider the Mackey-Glass equation, a delay differential equation for which it is still an open conjecture whether there exist parameters for which the system exhibits mathematically provable chaos. 
The phase space for this system is infinite dimensional, and we use a projection into $\R^8$ for our fractal dimension calculations. 
In \cite{farmer1982chaotic} Farmer reports the attractor's Lyapunov dimension to be $3.58 $, 
and we obtain a $\PH_0^{.5}$ dimension estimate close to this value. 
However, the spectrum of $\PH_0^\alpha$ dimension for various $\alpha$ had the largest spread of any of the fractals we studied, except the Rulkov attractor, ranging between $3.59$ for $ \alpha = .5$ and $4.06$ for $ \alpha = 3.5$ (see Figure \ref{fig:Attractor_spectrum}). 
While none of our dimension estimates appear to oscillate, all appear to converge slowly. The box-counting dimension estimation method described in~Section~\ref{sec:boxCompute} performed particularly poorly for this example, perhaps due to the high dimension and co-dimension. The estimate in Table~\ref{table:attractor_estimates} was computed by fitting a power law by hand. We did not attempt to compute either the $\PH_i$ dimensions for $1\leq i \leq 7,$ because the algorithm we use is impractical if the ambient dimension is greater than $3.$

\section{Earthquake Data}
\label{sec:earthquake}
We estimate the dimension of the Hauksson--Shearer Waveform Relocated Southern California earthquake catalog~\cite{2012hauksson,2007lin}. The hypocenter of an earthquake is the location beneath the earth's surface where the earthquake originates.  One can use the dimension estimate to study the geometry of earthquakes, i.e. by comparison to dimension estimates for fractures in rock surfaces or dislocations in crystals (see Section~6 of~\cite{2007kagan}). Previously, Harte~\cite{1998harte} estimated the correlation dimension of earthquake hypocenters in New Zealand and Japan, and Kagan~\cite{2007kagan} studied the correlation dimension of earthquakes in southern California and developed extensive heuristics to correct for errors in that computation.

\begin{figure}
\centering  
\subfigure[]{\includegraphics[width=.3\textwidth]{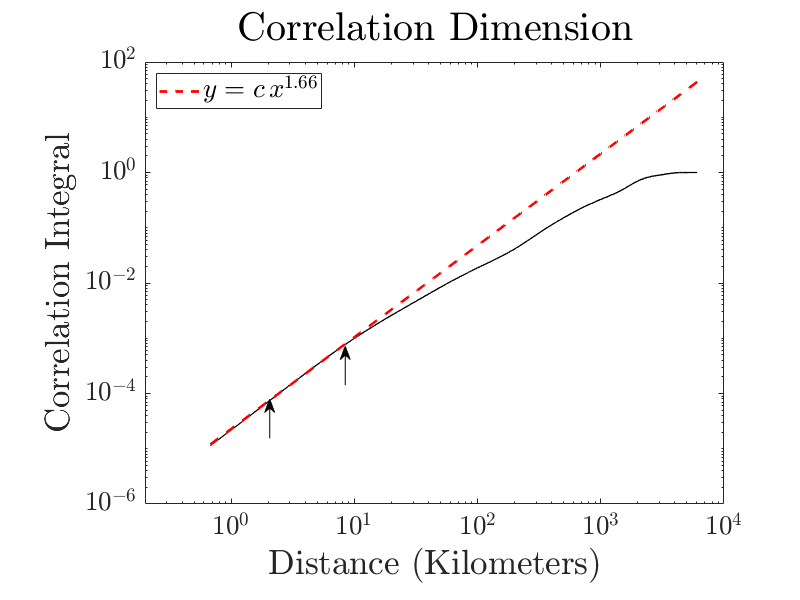}
\label{fig:earthquake_correlation_integral}}
\subfigure[]{
\includegraphics[width=.3\textwidth]{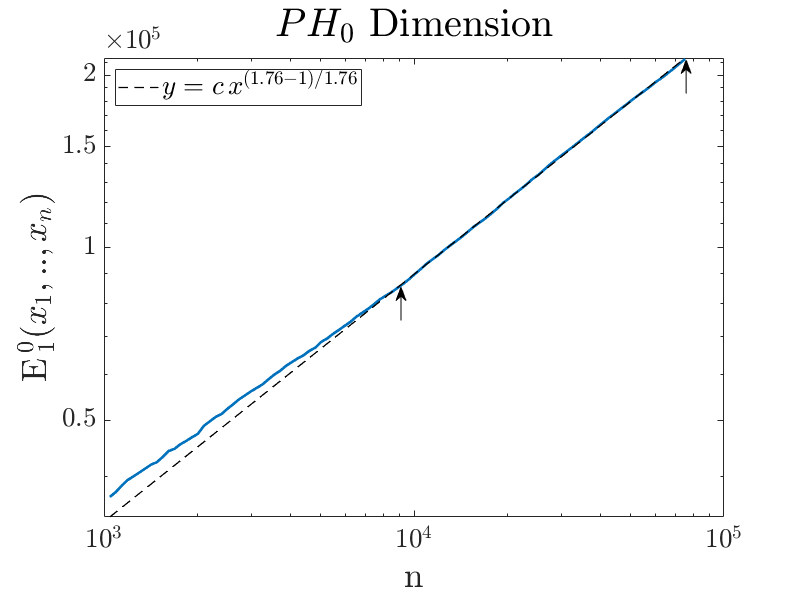}
\label{fig:earthquake_PHfit}}
\subfigure[]{\includegraphics[width=.3\textwidth]{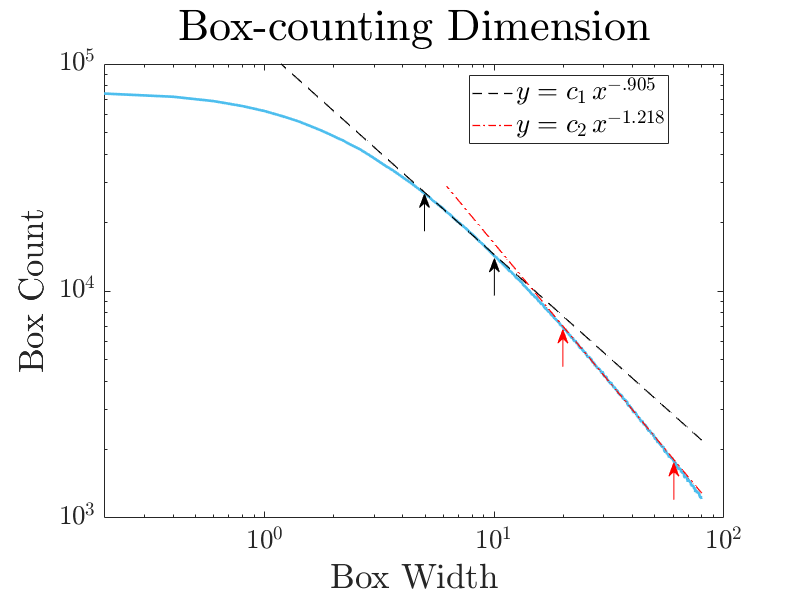}
\label{fig:earthquake_box}}
		\caption{\label{fig:earthquake}Dimension estimation for the earthquake hypocenter data set. Power laws were fitted in the ranges marked by the arrows.}
\end{figure}

We downloaded the coordinates of the hypocenters of the earthquakes in the catalogue from the Southern California Earthquake Data Center Website~\cite{2013scedc}, and selected the waveform-relocated earthquakes of magnitude greater than 2.0. This gave a data set of 74,929 earthquakes. The deepest earthquake in the data set is $33$ kilometers below sea level and the depth distribution has a single peak, so we did not split the data set into multiple samples as Harte did. The data set had small relocation errors, with the error in distance between nearby hypocenters estimated to be less than $.1$ kilometer for $90\%$ of such distances~\cite{2012hauksson}.

The presence of location errors means that a different methodology is required to estimate the correlation dimension than from the previously considered examples, where estimates at shorter lengthscales produce better results. To estimate the correlation dimension we plotted the correlation integral (Figure~\ref{fig:earthquake_correlation_integral}) and determined that the best range to fit a power law was between $2$ and $8.5$ kilometers. This yields a correlation estimate of about $1.66.$ This can be compared to Kagan's estimate of 1.5 for earthquakes from an earlier version of the same earthquake catalog that was processed with different methodology (Figure~8 of~\cite{2007kagan}), and to Harte's estimates of 1.8 and 1.5 for shallow earthquakes in Kanto, Japan and Wellington, New Zealand, respectively.

Figure~\ref{fig:earthquake_PHfit} shows the plot of $n$ versus $E_0^1\paren{x_1,\ldots,x_n}$ used to estimate the $\PH_0$ dimension. We estimated that the $\PH_0^1$ dimension is approximately $1.76$ when $\alpha=1.$ Similar computations for $\alpha=.5$ and $\alpha=1.5$ yielded dimension estimates of $1.75$ and $1.83,$ respectively. From this, and the comparison with the correlation dimension, we have evidence that the earthquake probability distribution is not regular. Also, note that  the $\PH_0^1$ dimension plot (Figure~\ref{fig:earthquake_PHfit}) appears to follow a power law in a large range, and the dimension can be estimated without fiddling with example-specific parameters. This precludes difficulties that arose in Harte's analysis where the correlation dimension estimate was very sensitive to the scale at which it was measured(Figure~20 in~\cite{1998harte}). Of course, the correlation dimension has an advantage of interpretability at different lengthscales --- the $x$-axis in Figure~\ref{fig:earthquake_correlation_integral} is kilometers, while it is number of samples in Figure~\ref{fig:earthquake_PHfit}.

The box-counting, $\PH_1,$ and $\PH_2$ dimensions did not produce good results. In the former case, it is unclear whether the plot of box width versus box count (Figure~\ref{fig:earthquake_box}) follows a power law at any range, but certainly no range shorter than $5$ kilometers. Fitting between $5$ and $10$ kilometers yielded an estimate of $.91$, which is very different than the estimates computed with other methods. Fitting at larger length scales, between $20$ and $50$ kilometers produced an estimate of $1.22.$ It is not surprising that the $\PH_1$ and $\PH_2$ dimensions produced poor results, given the small sample size. 

We also estimated  that $\text{comp}_{\PH_i}(X)$ likely equals zero for the earthquake data, for $i=0,1,2.$

\section{Conclusion}

Overall the performance of  our $\PH_0$  dimension calculations is comparable to the correlation dimension, and gives dimension estimates that are often greater for non-regular sets. Both the $\PH_0$ and correlation dimensions were more reliable than either the box-counting, $\PH_1,$ or $\PH_2$ dimensions. The correlation dimension and the $\PH_0$ dimension are easy to implement and quick to compute. The $\PH_0$ dimension may be more user-friendly in the sense that a single choice of power law range worked well for all examples.  Of these two, we do not view one dimension estimation technique as ``better'' than the other. Rather, they are complementary methods that --- if the underlying set is not regular --- can provide complementary information.

 In many applications of persistent homology, small intervals are often ignored in favor of larger features and discarded as ``noise''. However, as we have demonstrated here, there is sometimes signal in that noise. Adding the $\PH$ dimension to the current suite of $\PH$-based techniques used in applications such as machine learning will offer an orthogonal data descriptor.

\bibliographystyle{plain}
\bibliography{Bibliography}

\begin{thebibliography}{10}

\bibitem{2019adams}
Henry Adams, Manuchehr Aminian, Elin Farnell, Michael Kirby, Joshua Mirth,
  Rachel Neville, Chris Peterson, Patrick Shipman, and Clayton Shonkwiler.
\newblock A fractal dimension for measures via persistent homology.
\newblock {\em Abel Symposia}, 2019.

\bibitem{2017adams}
Henry Adams, Sofya Chepushtanova, Tegan Emerson, Eric Hanson, Michael Kirby,
  Francis Motta, Rachel Neville, Chris Peterson, Patrick Shipman, and Lori
  Ziegelmeier.
\newblock Persistence images: A stable vector representation of persistent
  homology.
\newblock {\em Journal of Machine Learning Research}, 2017.

\bibitem{1992aldous}
David Aldous and J.~Michael Steele.
\newblock Asymptotics for {E}uclidean minimal spanning trees on random points.
\newblock {\em Probability Theory and Related Fields}, 1992.

\bibitem{1987arneodo}
A.~Arneodo, G.~Grasseau, and E.~Kostelich.
\newblock Fractal dimensions and $f(\alpha)$ spectrum of the {H}{\'e}non
  attractor.
\newblock {\em Physics Letters A}, 1987.

\bibitem{1987badii}
R.~Badii and A.~Politi.
\newblock Renyi dimensions from local expansion rates.
\newblock {\em Physical Review A}, 1987.

\bibitem{2000baish}
James~W. Baish and Rakesh~K. Jain.
\newblock Fractals and cancer.
\newblock {\em Perspectives in Cancer Research}, 2000.

\bibitem{2000barbara}
Daniel Barbar\'{a} and Ping Chen.
\newblock Using the fractal dimension to cluster datasets.
\newblock In {\em KDD '00 Proceedings of the sixth ACM SIGKDD international
  conference on Knowledge discovery and data mining}, 2000.

\bibitem{2008beffara}
Vincent Beffara.
\newblock The dimension of the {SLE} curves.
\newblock {\em Annals of Probability}, 2008.

\bibitem{1996borovkova}
S.~Borovkova, R.~Burton, and H.~Dehling.
\newblock Consistency of the {T}akens estimator for the correlation dimension.
\newblock {\em Annals of Applied Probability}, 1996.

\bibitem{1928bouligand}
M.~G. Bouligand.
\newblock Ensembles impropres et nombre dimensionnel.
\newblock {\em Bulletin des Sciences Math{\'e}matiques}, 1928.

\bibitem{2015bubenik}
Peter Bubenik.
\newblock Statistical topological data analysis using persistence landscapes.
\newblock {\em Journal of Machine Learning Research}, 2015.

\bibitem{2011cagliaria}
Francesca Cagliaria and Claudia Landibc.
\newblock Finiteness of rank invariants of multidimensional persistent homology
  groups.
\newblock {\em Applied Mathematics Letters}, 2011.

\bibitem{2015camastra}
Francesco Camastra and Antonino Staiano.
\newblock Intrinsic dimension estimation: advances and open problems.
\newblock {\em Information Sciences}, 2015.

\bibitem{2016chazal}
F.~Chazal, V.~de~Silva, M.~Glisse, and S.~Oudot.
\newblock {\em The Structure and Stability of Persistence Modules}.
\newblock Springer, 2016.

\bibitem{2007cohensteiner}
David Cohen-Steiner, Herbert Edelsbrunner, and John Harer.
\newblock Stability of persistence diagrams.
\newblock {\em Discrete \& Computational Geometry}, 37(1), 2007.

\bibitem{2018curtin}
R.R. Curtin, M.~Edel, M.~Lozhnikov, Y.~Mentekidis, S.~Ghaisas, and S.~Zhang.
\newblock mlpack 3: a fast, flexible machine learning library.
\newblock {\em Journal of Open Source Software}, 2018.

\bibitem{1999davies}
S.~Davies and P.~Hall.
\newblock Fractal analysis of surface roughness by using spatial data.
\newblock {\em Journal of the Royal Statistical Society Series B}, 1999.

\bibitem{2007de_silva}
Vin {De Silva} and Robert Ghrist.
\newblock Coverage in sensor networks via persistent homology.
\newblock {\em Algebraic and Geometric Topology}, 2007.

\bibitem{2018divol}
Vincent Divol and Wolfgang Polonik.
\newblock On the choice of weight functions for linear representations of
  persistence diagrams.
\newblock arXiv:1807.03678, July 2018.

\bibitem{2002edelsbrunner}
H.~Edelsbrunner, D.~Letscher, and A.~Zomorodian.
\newblock Topological persistence and simplificitaion.
\newblock {\em Discrete and Computational Geometry}, 2002.

\bibitem{2008edelsbrunner}
Herbert Edelsbrunner and John Harer.
\newblock Persistent homology --- a survey.
\newblock {\em Contemporary Mathematics}, 2008.

\bibitem{2013edelsbrunner}
Herbert Edeslbrunner and Dmitriy Morozov.
\newblock Persistent homology: Theory and practice.
\newblock Technical Report LBNL-6037E, Lawrence Berkeley National Laboratory,
  2013.

\bibitem{2003edgar}
Gerald~A. Edgar.
\newblock {\em Classics on Fractals}.
\newblock Studies in Nonlinearity. Westview Press, 2003.

\bibitem{1985takens}
Takens F.
\newblock {\em Dynamical Systems and Bifurcations}, chapter On the numerical
  determination of the dimension of an attractor.
\newblock Springer, 1985.

\bibitem{2014falconer}
Kenneth Falconer.
\newblock {\em Fractal Geometry: Mathematical Foundations and Applications}.
\newblock Wiley, 2014.

\bibitem{farmer1982chaotic}
J~Doyne Farmer.
\newblock Chaotic attractors of an infinite-dimensional dynamical system.
\newblock {\em Physica D: Nonlinear Phenomena}, 4(3):366--393, 1982.

\bibitem{frederickson1983liapunov}
Paul Frederickson, James~L Kaplan, Ellen~D Yorke, and James~A Yorke.
\newblock The {L}iapunov dimension of strange attractors.
\newblock {\em Journal of differential equations}, 49(2):185--207, 1983.

\bibitem{2000gabella}
Marco Gabella, Sabastiano Pavone, , and Giovanni Perona.
\newblock Errors in the estimate of the fractal correlation dimension of
  raindrop spatial distribution.
\newblock {\em Journal of Applied Meteorology}, 2000.

\bibitem{garland2016exploring}
Joshua Garland, Elizabeth Bradley, and James~D Meiss.
\newblock Exploring the topology of dynamical reconstructions.
\newblock {\em Physica D: Nonlinear Phenomena}, 334:49--59, 2016.

\bibitem{2008ghrist}
Robert Ghrist.
\newblock Barcodes: the persistent homology of data.
\newblock {\em Bulletin of the American Mathematical Society}, 2008.

\bibitem{2015giusti}
Chad Giusti, Eva Pastalkova, Carina Curto, and Vladimir Itskov.
\newblock Clique topology reveals intrinsic geometric structure in neural
  correlations.
\newblock {\em Proceedings of the National Academy of Sciences}, 2015.

\bibitem{2012gneiting}
Tilmann Gneiting, Hana Sevcikova, and Donald Percival.
\newblock Estimators of fractal dimension: Assessing the roughness of time
  series and spatial data.
\newblock {\em Statistical Science}, 2012.

\bibitem{1984grassberger}
Peter Grassberger.
\newblock Generalizations of the {H}ausdorff dimension of fractal measures.
\newblock {\em Physics Letters A}, 1984.

\bibitem{1983grassberger}
Peter Grassberger and Itamar Procaccia.
\newblock Measuring the strangeness of strange attractors.
\newblock {\em Physica D: Nonlinear Phenomena}, 1983.

\bibitem{grebogi1984strange}
Celso Grebogi, Edward Ott, Steven Pelikan, and James~A Yorke.
\newblock Strange attractors that are not chaotic.
\newblock {\em Physica D: Nonlinear Phenomena}, 13(1-2):261--268, 1984.

\bibitem{1982greenside}
H.~S. Greenside, A.~Wolf, J.~Swift, and T.~Pignataro.
\newblock Impracticality of a box-counting algorithm for calculating the
  dimensionality of strange attractors.
\newblock {\em Physical Review A}, 1982.

\bibitem{1991guzzo}
Luigi Guzzo, Angela Iovino, Guido Chincarini, Riccardo Giovanelli, and
  Martha~P. Haynes.
\newblock Scale-invariant clustering in the large-scale distribution of
  galaxies.
\newblock {\em Astrophysical Journal, Part 2 - Letters}, 1991.

\bibitem{hammel1985global}
SM~Hammel, CKRT Jones, and Jerome~V Moloney.
\newblock Global dynamical behavior of the optical field in a ring cavity.
\newblock {\em JOSA B}, 2(4):552--564, 1985.

\bibitem{1998harte}
D.~Harte.
\newblock Dimension estimates of earthquake epicentres and hypocentres.
\newblock {\em Journal of Nonlinear Science}, 1998.

\bibitem{2012hauksson}
Egill Hauksson, Wenzheng Yang, and Peter~M. Sheare.
\newblock Waveform relocated earthquake catalog for southern {C}alifornia (1981
  to june 2011)r.
\newblock {\em Bulletin of the Seismological Society of America}, 2012.

\bibitem{1918hausdorff}
Felix Hausdorff.
\newblock Dimension und {\"a}u{\ss}eres ma{\ss}.
\newblock {\em Mathematische Annalen}, 1918.

\bibitem{henon1976two}
Michel H{\'e}non.
\newblock A two-dimensional mapping with a strange attractor.
\newblock In {\em The Theory of Chaotic Attractors}, pages 94--102. Springer,
  1976.

\bibitem{2016hiraoka}
Yasuaki Hiraoka, Takenobu Nakamura, Akihiko Hirata, Emerson~G Escolar, Kaname
  Matsue, and Yasumasa Nishiura.
\newblock Hierarchical structures of amorphous solids characterized by
  persistent homology.
\newblock {\em Proceedings of the National Academy of Sciences}, 2016.

\bibitem{2003hu}
Song Hu, Xuexing Sun, Jun Xiang, Min Li, Peisheng Li, and Liqi Zhang.
\newblock Correlation characteristics and simulations of the fractal structure
  of coal char.
\newblock {\em Communications in Nonlinear Science and Numerical Simulation},
  2003.

\bibitem{ibarz2011map}
Borja Ibarz, Jos{\'e}~Manuel Casado, and Miguel~AF Sanju{\'a}n.
\newblock Map-based models in neuronal dynamics.
\newblock {\em Physics reports}, 501(1-2):1--74, 2011.

\bibitem{jorba2008mechanism}
{\`A}ngel Jorba and Joan~Carles Tatjer.
\newblock A mechanism for the fractalization of invariant curves in
  quasi-periodically forced 1-d maps.
\newblock {\em Discrete \& Continuous Dynamical Systems-B}, 10:537--567, 2008.

\bibitem{2007kagan}
Yan~Y. Kagan.
\newblock Earthquake spatial distribution: the correlation dimension.
\newblock {\em Geophysical Journal International}, 2007.

\bibitem{kaplan1979chaotic}
James~L Kaplan and James~A Yorke.
\newblock Chaotic behavior of multidimensional difference equations.
\newblock In {\em Functional differential equations and approximation of fixed
  points}, pages 204--227. Springer, 1979.

\bibitem{1996kesten}
Harry Kesten and Sungchul Lee.
\newblock The central limit theorem for weighted minimal spanning trees on
  random points.
\newblock {\em Annals of Applied Probability}, 1996.

\bibitem{1989leibovitch}
Larry~S. Liebovitch and Tibor Toth.
\newblock A fast algorithm to determine fractal dimensions by box counting.
\newblock {\em Physics Letters A}, 1989.

\bibitem{2007lin}
G.~Lin, P.~M. Shearer, and E.~Hauksson.
\newblock Applying a three-dimensional velocity model, waveform cross
  correlation, and cluster analysis to locate southern {C}alifornia seismicity
  from 1981 to 2005.
\newblock {\em Journal of Geophysical Research: Solid Earth}, 2007.

\bibitem{2009lopes}
R.~Lopes and N.~Betrouni.
\newblock Fractal and multifractal analysis: A review.
\newblock {\em Medical Image Analysis}, 2009.

\bibitem{lorenz1963deterministic}
Edward~N Lorenz.
\newblock Deterministic nonperiodic flow.
\newblock {\em Journal of the atmospheric sciences}, 20(2):130--141, 1963.

\bibitem{1990lovejoy}
Shaun Lovejoy and Daniel Schertzer.
\newblock Fractals, raindrops and resolution dependence of rain measurements.
\newblock {\em Journal of Applied Meteorology}, 1990.

\bibitem{mackey1977oscillation}
Michael~C Mackey and Leon Glass.
\newblock Oscillation and chaos in physiological control systems.
\newblock {\em Science}, 197(4300):287--289, 1977.

\bibitem{2012macpherson}
R.~D. MacPherson and B.~Schweinhart.
\newblock Measuring shape with topology.
\newblock {\em Journal of Mathematical Physics}, 53(7), 2012.

\bibitem{1977mandelbrot}
Beno\^{i}t Mandelbrot.
\newblock {\em Fractals: Form, Chance and Dimension}.
\newblock W.H.Freeman \& Company, 1977.

\bibitem{1982mandelbrot}
Beno\^{i}t Mandelbrot.
\newblock {\em The Fractal Geometry of Nature}.
\newblock W. H. Freeman and Co., 1982.

\bibitem{2010march}
William~B. March, Parikshit Ram, and Alexander~G. Gray.
\newblock Fast {E}uclidean minimum spanning tree: Algorithm, analysis, and
  applications.
\newblock {\em KDD '10 Proceedings of the 16th ACM SIGKDD international
  conference on Knowledge discovery and data mining}, 2010.

\bibitem{GUDHI_persistence}
Cl\'ement Maria.
\newblock Persistent cohomology.
\newblock In {\em {GUDHI} User and Reference Manual}. {GUDHI Editorial Board},
  2015.

\bibitem{1993martinez}
Vicent~J. Martinez, R.~Dominguez-Tenreiro, and L.~J. Roy.
\newblock {H}ausdorff dimension from the minimal spanning tree.
\newblock {\em Physical Review E}, 1992.

\bibitem{mischaikow1999construction}
Konstantin Mischaikow, Marian Mrozek, J~Reiss, and Andrzej Szymczak.
\newblock Construction of symbolic dynamics from experimental time series.
\newblock {\em Physical Review Letters}, 82(6):1144, 1999.

\bibitem{2012mo}
Dengyao Mo and Samuel~H. Huang.
\newblock Fractal-based intrinsic dimension estimation and its application in
  dimensionality reduction.
\newblock {\em IEEE Transactions on Knowledge and Data Engineering}, 2012.

\bibitem{myers2019persistent}
Audun Myers, Elizabeth Munch, and Firas~A Khasawneh.
\newblock Persistent homology of complex networks for dynamic state detection.
\newblock {\em arXiv preprint arXiv:1904.07403}, 2019.

\bibitem{1990nerenberg}
M.~A.~H. Nerenberg and Christopher Essex.
\newblock Correlation dimension and systematic geometric effects.
\newblock {\em Physical Review A}, 1990.

\bibitem{2013peng}
Xin Peng, Wei Qi, Mengfan Wang, Rongxin Su, and Zhimin He.
\newblock Backbone fractal dimension and fractal hybrid orbital of protein
  structure.
\newblock {\em Communications in Nonlinear Science and Numerical Simulation},
  2013.

\bibitem{1970renyi}
Alfred Renyi.
\newblock {\em Probability Theory}.
\newblock North Holland Publishing Company, 1970.

\bibitem{1991rieu}
Michel Rieu and Garrison Sposito.
\newblock Fractal fragmentation, soil porosity, and soil water properties: I.
  theory.
\newblock {\em Soil Science Society of America Journal}, 1991.

\bibitem{2000robins}
Vanessa Robins.
\newblock {\em Computational topology at multiple resolutions: foundations and
  applications to fractals and dynamics}.
\newblock PhD thesis, University of Colorado at Boulder, 2000.

\bibitem{2018rosenberg}
Eric Rosenberg.
\newblock {\em A Survey of Fractal Dimensions of Networks}.
\newblock Springer, 2018.

\bibitem{GUDHI_alpha}
Vincent Rouvreau.
\newblock Alpha complex.
\newblock In {\em {GUDHI} User and Reference Manual}. {GUDHI Editorial Board},
  2015.

\bibitem{rulkov2001regularization}
Nikolai~F Rulkov.
\newblock Regularization of synchronized chaotic bursts.
\newblock {\em Physical Review Letters}, 86(1):183, 2001.

\bibitem{2017saadatfar}
Mohammad Saadatfar, Hiroshi Takeuchi, Vanessa Robins, Nicolas Francois, and
  Yasuaki Hiraoka.
\newblock Pore configuration landscape of granular crystallization.
\newblock {\em Nature Communications}, 8, 2017.

\bibitem{1994sarkar}
N.~Sarkar and B.B. Chaudhuri.
\newblock An efficient differential box-counting approach to compute fractal
  dimension of image.
\newblock {\em IEEE Transactions on Systems, Man, and Cybernetics}, 1994.

\bibitem{sauer1991embedology}
Tim Sauer, James~A Yorke, and Martin Casdagli.
\newblock Embedology.
\newblock {\em Journal of statistical Physics}, 65(3-4):579--616, 1991.

\bibitem{2013scedc}
SCEDC.
\newblock Southern california earthquake center., 2013.
\newblock Caltech.Dataset.

\bibitem{2019schweinhart}
B.~Schweinhart.
\newblock The persistent homology of random geometric complexes on fractals.
\newblock arXiv:1808.02196.

\bibitem{2018schweinhart}
Benjamin Schweinhart.
\newblock Persistent homology and the upper box dimension.
\newblock arXiv:1802.00533, 2018.

\bibitem{sprott2001improved}
Julien~Clinton Sprott and George Rowlands.
\newblock Improved correlation dimension calculation.
\newblock {\em International Journal of Bifurcation and Chaos},
  11(07):1865--1880, 2001.

\bibitem{1988steele}
J.~Michael Steele.
\newblock Growth rates of {E}uclidean minimal spanning trees with power
  weighted edges.
\newblock {\em Annals of Probability}, 1988.

\bibitem{1980takens}
Floris Takens.
\newblock {\em Dynamical Systems and Turbulence}, chapter Detecting Strange
  Attractors in Turbulence, pages 366--381.
\newblock Springer, 1980.

\bibitem{1991taylor}
Charles~C. Taylor and James Taylor.
\newblock Estimating the dimension of a fractal.
\newblock {\em Journal of the Royal Statistical Society. Series B
  (Methodological)}, 1991.

\bibitem{1990theiler}
James Theiler.
\newblock Estimating fractal dimension.
\newblock {\em Journal of the Optical Society of America A}, 1990.

\bibitem{1990theiler_b}
James Theiler.
\newblock Statistical precision of dimension estimators.
\newblock {\em Physical Review A}, 1990.

\bibitem{theiler1988quantifying}
James~Patrick Theiler.
\newblock {\em Quantifying Chaos: Practical Estimation of the Correlation
  Dimension.}
\newblock PhD thesis, California Institute of Technology, 1988.

\bibitem{2010traina}
Caetano Traina, Agma Traina, Leejay Wu, and Christos Faloutsos.
\newblock Fast feature selection using fractal dimension.
\newblock {\em Journal of Information and Data Management}, 2010.

\bibitem{1992weygaert}
Rien van~de Weygaert, Bernard~J.T. Jones, and Vincent~J. Martinez.
\newblock The minimal spanning tree as an estimator for generalized dimensions.
\newblock {\em Physics Letters A}, 1992.

\bibitem{1927vietoris}
Leopold Vietoris.
\newblock {\"U}ber den h{\"o}heren zusammenhang kompakter r{\"a}ume und eine
  klasse von zusammenhangstreuen abbildungen.
\newblock {\em Mathematische Annalen}, 1927.

\bibitem{2014xia}
Kelin Xia and Guo-Wei Wei.
\newblock Persistent homology analysis of protein structure, flexibility, and
  folding.
\newblock {\em International journal for numerical methods in biomedical
  engineering}, 2014.

\bibitem{young1982dimension}
Lai-Sang Young.
\newblock Dimension, entropy and {L}yapunov exponents.
\newblock {\em Ergodic theory and dynamical systems}, 2(1):109--124, 1982.

\bibitem{young2013mathematical}
Lai-Sang Young.
\newblock Mathematical theory of {L}yapunov exponents.
\newblock {\em Journal of Physics A: Mathematical and Theoretical},
  46(25):254001, 2013.

\bibitem{2008yu}
Boming Yu.
\newblock Analysis of flow in fractal porous media.
\newblock {\em Applied Mechanics Reviews}, 2008.

\bibitem{2000yukich}
J.E. Yukich.
\newblock Asymptotics for weighted minimal spanning trees on random points.
\newblock {\em Stochastic Processes and their Applications}, 2000.

\end{thebibliography}

\appendix

\section{Self-Similar Fractals}
\label{appendix:self_similar}
\subsection{The Sierpinski Triangle}
The Sierpinski triangle (Figure~\ref{fig:sierpinski}) is defined by iteratively removing equilateral triangles from a larger equilateral triangle. Begin by sub-dividing the equilateral triangle formed by $(0,0),(1,0),$ and $(\sqrt{3}/2,1/2)$ into four congruent triangles, and removing the center triangle. Repeat this process on the remaining three triangles, and continue ad infinitum. The Sierpinski triangle equals three copies of itself rescaled by the factor $1/2$ so the dimension of the resulting set is $\frac{\log\paren{3}}{\log\paren{2}}$ (for a self-similar set equal to $m$ copies of itself rescaled by a factor $r,$ the self-similarity dimension equals $\frac{\log\paren{m}}{\log\paren{1/r}}.$

We sampled points from the natural measure on the Sierpinski Triangle by sampling random integers $\set{a_1, a_2, \ldots}\in\set{0,1,2}$ and computing
\[\paren{x,y}=\paren{\sum_{i:a_i=1}2^{-i}+ \sum_{i:a_i=2}2^{-i-1}, \sum_{i:a_i=2}\sqrt{3}\,2^{-i-1}}\,.\]
In practice, we end this procedure at $i=64.$ 

\subsection{The Cantor Dust}
The standard middle-thirds Cantor set $C$ is the set formed by removing the interval $\paren{1/3,2/3}$ from the closed interval $\brac{0,1},$ and iteratively removing the middle third of the remaining intervals. The dimension of $C$ is $\frac{\log\paren{2}}{\log\paren{3}}.$ The Cantor dust (Figure~\ref{fig:cantordust}) is the product $C\times C;$ its dimension is  $\frac{2\log\paren{2}}{\log\paren{3}}.$

We sampled points from the natural measure on $C$ by sampling random integers $\set{a_1, a_2, \ldots}\in\set{0,1}$ and computing
\[x=2\sum_{i=1}^\infty a_i 3^{-i}\,.\]
In practice, we truncated the summation at $i=64.$ A point $(x_1,x_2)$ from the natural measure on the Cantor dust can then be sampled by independently sampling $x_1$ and $x_2$ as above.

\subsection{The Cantor Set Cross an Interval}

Consider the set $C\times \brac{0,1},$ where $C$ is the Cantor set defined in the previous section. The dimension of this set is $1+\frac{\log\paren{2}}{\log\paren{3}},$ one greater than the dimension of the Cantor set.

We sampled points $(x,y)$ from the natural measure on the Cantor set cross an interval by sampling a random point $x$ on the Cantor set by the procedure described in the previous section, and a random real number $y$ from the uniform distribution on $\brac{0,1}.$

\subsection{The Menger Sponge}
The Menger sponge is defined by iteratively removing cubes from a larger cube. Begin by sub-dividing the unit cube in $\mathbb{R}^3$ into nine smaller sub-cubes, and remove nine of them: one from the center of the original cube one from the center of each of the eight faces. Repeat this process on each of the $20$ remaining sub-cubes and continue ad infinitum. The resulting structure is equal to $20$ copies of itself rescaled by a factor $1/3$ so the dimension of the  resulting set is $\frac{\log\paren{3}}{\log\paren{20}}.$

For $j\in\mathbb{n}$ we sample a uniform three-tuple of integers $\paren{x_i,y_i,z_i}\in\set{0,1,2}$ throwing out and re-sampling any three-tuple for which two or more of the coordinates equals one. To sample a point from the natural measure on the Menger sponge, we form the sum
\[\paren{x,y,z}=\sum_{i=0}^\infty 3^{-i}\paren{x_i ,y_i ,z_i}\,.\]

In practice, we end this procedure at $i=64.$

\section{Chaotic Attractors}
\label{sec:Attractors_Description}

\subsection{H\'enon Map}
\label{sec:Henon}
The H\'enon map \cite{henon1976two} is given by $(x_n,y_n) \mapsto (x_{n+1}, y_{n+1})$ where: 
\begin{align*}
x_{n+1} &=  1 - a x_n^2 + y_n \\
y_{n+1} &= b  x_n .
\end{align*}
We used the  parameters of $ a = 1.4$ and $ b = 0.3$.  
We generated $50$ time series using randomly chosen initial conditions, computing a trajectory of length $1.1 \cdot 10^{7}$ and discarding the initial $10^6$ points in the series.

\subsection{Ikeda Map}
\label{sec:Ikeda}
The complex Ikeda map \cite{hammel1985global} is given by  $z_n \mapsto z_{n+1}$ where:
\begin{align*}
z_{n+1} = a + R  \ \mbox{exp}\left[ i \left( \phi -\frac{p}{1+|z_{n}|^2} \right) \right] z_{n}. 
\end{align*}
We used the parameters: 
\begin{align*}
a &= 1 &
R &= 0.9 &
\phi &= 0.4 &
p  &= 6. 
\end{align*}
We generated 50 time series using  randomly chosen initial conditions, computing a trajectory of length $1.1 \cdot 10^{6}$ and discarding the initial $10^5$ points in the series.

\subsection{Rulkov Map}
\label{sec:Rulkov}
The chaotic Rulkov map \cite{rulkov2001regularization,ibarz2011map} is given by $ (x_n , y_n) \mapsto (x_{n+1} , y_{n+1})$ where:  
\begin{align*}
x_{n+1} &= \frac{\alpha }{1 + x_n^2} + y_n \\
y_{n+1} &= y_n - \mu(x_n-\sigma)
\end{align*}
We use the parameters  
\begin{align*}
\mu &= 0.0001 &
\alpha &= 3.75 &
\sigma &= -1 &
I &= 0
\end{align*}
We generated 50 time series using randomly chosen initial conditions, computing a trajectory of length $1.1 \cdot 10^{6}$ and discarding the initial $10^5$ points in the series.

\subsection{Lorenz System}
\label{sec:Lorenz}
The Lorenz system \cite{lorenz1963deterministic}  given by the system of ordinary differential equations below:
\begin{align*}
\dot{x} &= \sigma(y-x)  \\
\dot{y} &= x(\rho - z) - y\\
\dot{z} &= xy - \beta z
\end{align*}
We use the parameters: 
\begin{align*}
\rho &= 28 &  \sigma & = 10 & \beta &= 8/3.
\end{align*} 

We generated 50 time series using randomly chosen initial conditions, integrating forward for $1.01 \times 10^5$  units of time,  and discarding the initial $10^3$ units of time.  
We used MATLAB's \texttt{ode45}  integrator with a relative error  tolerance of $10^{-9}$ and an absolute error tolerance of $10^{-9}$. 
We sampled our trajectories at a rate of $10$Hz, producing time series with $10^6$ points each.

\subsection{Mackey-Glass} 
\label{sec:MackeyGlass}
The Mackey-Glass equation \cite{mackey1977oscillation} is given by the scalar delay differential equation below:
\begin{align}
y'(t)  &= -a y(t) + b
\frac{y(t - \tau)}{  1+ y(t-\tau)^n}.
\end{align}
We take parameters:  
\begin{align*}
a&= 1 &
b&= 2 &
\tau&= 3 &
n&=10 .
\end{align*}
A natural phase space for this dynamical system is the Banach space of continuous functions $C =  C^0( [-\tau,0] , \R)$. 
For our dimension calculations we discretize this space using a projection map $\pi_m :C \to \R^{m}$ (for $m=8$) which evaluates a function $y \in C$ at $m$ maximally spaced points on the interval $ [-\tau , 0]$.

We generated 25 time series using  randomly chosen initial conditions and integrating forward for $ \tfrac{\tau}{m-1} \cdot 10^6 + 10^4 $ units of time.
We used MATLAB's \texttt{dde23}  integrator with the default error tolerances for the first $5,000$ units of time, and the remainders of each trajectory were computed using a relative error tolerance of $10^{-5}$ and an absolute error tolerance of $10^{-9}$. 
The initial transient periods of $10^4$ units of time were discarded, and the remainders were sampled at a rate of $\tfrac{\tau}{m-1}$Hz, producing time series with $10^6$ points each.

\end{document}